\documentclass[10pt]{amsart}
\usepackage{amssymb, amsthm, amsmath, amscd}

\theoremstyle{definition}
\newtheorem{thm}{Theorem}[section]
\newtheorem{lem}[thm]{Lemma}
\newtheorem{prp}[thm]{Proposition}
\newtheorem{dfn}[thm]{Definition}
\newtheorem{cor}[thm]{Corollary}

\newtheorem{rmk}[thm]{Remark}
\newtheorem{ntn}[thm]{Notation}
\newtheorem{exa}[thm]{Example}

\newcommand{\af}{\alpha}
\newcommand{\bt}{\beta}
\newcommand{\gm}{\gamma}
\newcommand{\dt}{\delta}
\newcommand{\ep}{\varepsilon}
\newcommand{\zt}{\zeta}
\newcommand{\et}{\eta}
\newcommand{\ch}{\chi}
\newcommand{\io}{\iota}
\newcommand{\te}{\theta}
\newcommand{\ld}{\lambda}
\newcommand{\sm}{\sigma}

\newcommand{\ph}{\varphi}

\newcommand{\rh}{\rho}
\newcommand{\om}{\omega}
\newcommand{\ta}{\tau}

\newcommand{\Gm}{\Gamma}
\newcommand{\Dt}{\Delta}

\newcommand{\rsz}[1]{\raisebox{0ex}[0.8ex][0.8ex]{$#1$}}
\renewcommand{\phi}{\varphi}

\newcommand{\N}{\mathbb{N}}
\newcommand{\Z}{\mathbb{Z}}
\newcommand{\Q}{\mathbb{Q}}
\newcommand{\R}{\mathbb{R}}
\newcommand{\C}{\mathbb{C}}

\newcommand{\Hom}{\operatorname{Hom}}

\newcommand{\Aff}{\operatorname{Aff}}

\newcommand{\Tr}{\operatorname{Tr}}
\newcommand{\id}{\operatorname{id}}

\newcommand{\dist}{\operatorname{dist}}

\newcommand{\sa}{\mathrm{sa}}

\newcommand{\ata}{AT~algebra}
\newcommand{\aha}{AH~algebra}

\newcommand{\ca}{C*-algebra}
\newcommand{\csa}{C*-subalgebra}
\newcommand{\uct}{Universal Coefficient Theorem}
\newcommand{\ct}{continuous}
\newcommand{\cfn}{continuous function}
\newcommand{\pj}{projection}

\newcommand{\cms}{compact metric space}

\newcommand{\nbhd}{neighborhood}

\newcommand{\hm}{homomorphism}

\newcommand{\ifo}{if and only if}
\newcommand{\rsha}{recursive subhomogeneous algebra}

\newcommand{\rshd}{recursive subhomogeneous decomposition}
\newcommand{\mops}{mutually orthogonal \pj s}

\newcommand{\mh}{minimal homeomorphism}
\newcommand{\hme}{homeomorphism}

\newcommand{\tgca}{transformation group \ca}
\newcommand{\Tgca}{Transformation group \ca}

\newcommand{\sint}{{\mathrm{int}}}
\newcommand{\diag}{{\mathrm{diag}}}
\newcommand{\diam}{{\mathrm{diam}}}
\newcommand{\supp}{{\mathrm{supp}}}

\newcommand{\ts}[1]{{\textstyle{#1}}}
\newcommand{\ds}[1]{{\displaystyle{#1}}}
\newcommand{\dirlim}{\varinjlim}

\newcommand{\andeqn}{\,\,\,\,\,\, {\mbox{and}} \,\,\,\,\,\,}

\renewcommand{\S}{\subset}
\newcommand{\SM}{\setminus}
\newcommand{\E}{\varnothing}
\newcommand{\I}{\infty}

\title{Crossed products by minimal homeomorphisms}
\author{Huaxin Lin}
\author{N.\  Christopher Phillips}
\thanks{Research of the first author
partially supported by NSF grants DMS 009703 and DMS 0355273,
and the Shanghai Priority Academic Disciplines.
Research of the second author
partially supported by NSF grant DMS 0302401.}
\date{14 Aug.~2004}

\begin{document}

\begin{abstract}
Let $X$ be an infinite compact metric space with finite covering dimension
and let $h \colon X \to X$ be a minimal homeomorphism.
We show that
the associated crossed product \ca\   $A = C^* (\Z, X, h)$
has tracial rank zero
whenever the image of $K_0 (A)$ in $\Aff (T (A))$ is dense.
As a consequence, we show that these crossed product \ca s are
in fact simple \aha s with real rank zero.
When $X$ is connected and $h$ is further assumed to
be uniquely ergodic,
then the above happens if and only if the rotation number associated
to $h$ has irrational values.

By applying the classification result for
nuclear simple \ca s with tracial rank zero, we show that two
such dynamical systems have isomorphic crossed products
if and only if they have isomorphic scaled ordered $K$-theory.
\end{abstract}

\maketitle

\section{Introduction}

\Tgca s of \mh s of compact metric spaces have a
long history in the theory of \ca s.
The irrational rotation algebras,
among the most prominent examples in the whole subject,
have this form.
The remarkable work of Giordano, Putnam, and Skau~\cite{GPS}
shows that the \tgca s of two \mh s of the Cantor set are
isomorphic \ifo\  the \hme s are strong orbit equivalent.
The \tgca s of Furstenberg transformations on toruses have also
attracted considerable attention; as just a few examples,
we mention~\cite{Pc86}, \cite{Ji}, \cite{Kd}, and~\cite{RM}.
Connes' example (in Section~5 of~\cite{Cn0})
of a simple \ca\  with no nontrivial \pj s
is the \tgca\  of a minimal diffeomorphism of $S^3.$

One naturally wants to understand the structure of such \ca s.
That the \tgca s of \mh s of the Cantor set are
\ata s (direct limits of circle algebras)
with real rank zero is implicit in Section~8 of~\cite{HPS},
with the main step having been done in~\cite{Pt2}.
Elliot and Evans proved in~\cite{EE} that the irrational
rotation algebras are \ata s with real rank zero.
In particular, these algebras all belong to the class
known currently to be classifiable by K-theoretic invariants
in the sense of the Elliott classification program~\cite{El5}.
It is proved in a very long and still unpublished paper~\cite{LP3}
(see the survey article~\cite{LP2})
that the \tgca\  of a minimal diffeomorphism
$h$ of a compact
smooth manifold $X$ is a direct limit, with no dimension growth,
of recursive subhomogeneous \ca s.
Using this result, one can deduce that if the K-theory and
traces on $C^* (\Z, X, h)$ are consistent with this algebra
having real rank zero (in the sense described in the next paragraph),
then $A$ has tracial rank zero (is tracially AF)
in the sense of \cite{Ln1},~\cite{Ln2}
(see Definition~\ref{DTR0} below),
and is hence classifiable.

In this paper, we prove the following.
Let $X$ be an infinite compact metric space
with finite covering dimension,
and let $h \colon X \to X$ be a \mh.
Let $A = C^* (\Z, X, h)$ be the \tgca.
Suppose (notation explained below) that the map
$\rh_A \colon K_0 (A) \to \Aff (T (A))$
has dense range.
(It is well known, and proved in Proposition~1.10(a) of~\cite{Ph8},
that this condition is necessary for $A$ to have real rank zero.)
Then in fact $A$ has tracial rank zero.
Here $T (A)$ is the space of tracial states on $A$
with the weak* topology,
and if $\Dt$ is a compact convex set then
$\Aff (\Dt)$ is the space of real valued affine \cfn s on $\Dt.$
The map $\rh_A \colon K_0 (A) \to \Aff (T (A))$ is defined by
$\rh_A (\et) (\ta) = \ta_* (\et)$ for $\et \in K_0 (A)$ and
$\ta \in T (A).$

It follows that $A$ is classifiable,
and is in fact a simple unital \aha\  with real rank zero.
Furthermore, if the \tgca s of two \mh s of finite dimensional
compact metric spaces, both satisfying the dense range condition,
have isomorphic scaled ordered K-theory, then they are isomorphic.
This result gives new examples of distinct \mh s with
isomorphic \tgca s.

The proof has the advantage over~\cite{LP3} of being short,
and of not requiring any smoothness on $X$ or on $h.$
It has the disadvantage of giving no information about
the \tgca\  when $\rh_{C^* (\Z, X, h)}$ does not have dense range.
The direct limit decomposition in~\cite{LP3} requires no such
assumption,
and implies that the \tgca\  has stable rank one and that the
order on \pj s is determined by traces even when it does not have
real rank zero,
such as for the example of Section~5 of~\cite{Cn0}.

To describe the main idea of the proof,
let $u \in C^* (\Z, X, h)$ be the canonical unitary,
and let $y \in X.$
Consider the subalgebra
$B = C^* (\Z, X, h)_{ \{ y \} }$ of $C^* (\Z, X, h)$
generated by $C (X)$ and all $u f$ with $f \in C (X)$ and $f (y) = 0.$
(This subalgebra plays a key role in~\cite{LP1},
and is a generalization of a subalgebra
originally introduced in~\cite{Pt1}.)
Under our hypotheses,
$\rh_B \colon K_0 (B) \to \Aff (T (B))$ also has dense range.
If there are only countably many ergodic $h$-invariant measures,
then Theorem~4.4 of~\cite{Ph8}
implies that $B$ has tracial rank zero.
Even if this is not the case, one can combine intermediate
results from~\cite{Ph8} and dynamical arguments to show that
$B$ has tracial rank zero for ``most'' (presumably actually all)
choices of $y.$
Having this, one can use Berg's technique to find arbitrarily ``small''
approximately $h$-invariant (in a suitable sense) \pj s $p$
in $C^* (\Z, X, h)$ and functions $f \in C (X)$ such that
$f = 1$ on a \nbhd\  of $y$ and such that $p f = f.$
The corners $(1 - p) C^* (\Z, X, h) (1 - p)$ are,
in a suitable sense, approximately contained in $B.$
Since $B$ has tracial rank zero and $1 - p$ is ``large'',
one can then deduce that $C^* (\Z, X, h)$ has tracial rank zero.

This paper is organized as follows.
In Section~\ref{Sec:SA} we introduce notation and give some
elementary properties of algebras related to the subalgebra $B$
in the description above.
In Section~\ref{Sec:TRSA}, we prove that,
under our hypotheses, $B$ (usually) has tracial rank zero.
Section~\ref{Sec:TRCrP} contains the proof that $C^* (\Z, X, h)$
has tracial rank zero.
In Section~\ref{Sec:Examples},
we use our result to examine some examples.
In particular, we prove isomorphism of the \tgca s in some
pairs of examples in~\cite{Ph} for which the isomorphism
question was left open there, due to lack of smoothness.
Section~\ref{Sec:AppConj} considers approximate conjugacy of the
\hme s in some of these examples.

We use the notation $p \precsim q$ to mean that a \pj\  $p$
is Murray-von Neumann equivalent to a sub\pj\  of a \pj\  $q.$
By convention,
if $B$ is said to be a unital subalgebra of a unital \ca\  $A,$
then, unless otherwise specified,
we mean that the identity of $B$ is the same as that of $A.$

For the convenience of the reader,
we recall the meaning of tracial rank zero
(or tracial topological rank zero) for simple \ca s.

\begin{dfn}\label{DTR0}
Let $A$ be a simple unital \ca.
Then $A$ has tracial rank zero if
for every finite subset ${\mathcal{F}} \S A,$ every $\ep > 0,$
and every nonzero positive element $c \in A,$
there exists a \pj\  $p \in A$ and a unital finite dimensional
subalgebra $E \S p A p$ such that:
\begin{itemize}
\item[(1)]
$\| [a, p] \| < \ep$ for all $a \in {\mathcal{F}}.$
\item[(2)]
$\dist (p a p, E) < \ep$ for all $a \in {\mathcal{F}}.$
\item[(3)]
$p$ is Murray-von Neumann equivalent to a \pj\  in ${\overline{c A c}}.$
\end{itemize}
\end{dfn}

That this definition is equivalent to the original one follows
from Proposition~3.8 of~\cite{Ln1}.

\vspace{1ex}

{\textbf{Acknowledgements:}}
The second author would like to thank
Mike Boyle and Tomasz Downarowicz
for valuable email correspondence
on \mh s of finite dimensional compact metric spaces.
In particular, Mike Boyle called attention to the reference~\cite{Kl}.

\section{Subalgebras of the transformation group
                            C*-algebra}\label{Sec:SA}

Let $X$ be a compact metric space,
and let $h \colon X \to X$ be a homeomorphism.
Then the induced automorphism $\af$ of $C (X)$ is
$\af (f) = f \circ h^{-1}.$
In the \tgca\  $C^* (\Z, X, h),$ we normally write
$u$ for the standard unitary representing the generator of $\Z.$
Then $u f u^* = \af (f) = f \circ h^{-1}.$

\begin{ntn}\label{SubalgNtn}
Let $X$ be a compact metric space,
and let $h \colon X \to X$ be a homeomorphism.
For a nonempty closed subset $Y \subset X,$ we define
the C*-subalgebra $C^* (\Z, X, h)_Y$ to be
\[
C^* (\Z, X, h)_Y
     = C^* ( C (X), u C_0 (X \setminus Y)) \subset C^* (\Z, X, h).
\]
We will often let $A$ denote the
\tgca\  $C^* (\Z, X, h),$ in which case we refer to $A_Y.$
\end{ntn}

The case of direct use to us is $Y = \{ y \}.$
The following immediate structural fact will be crucial.

\begin{rmk}\label{LimSubalg}
Let the notation be as in Notation~\ref{SubalgNtn},
with $A = C^* (\Z, X, h).$
If $Y_1 \supset Y_2 \supset \cdots$ is a decreasing sequence
of closed subsets of $X$ with $\bigcap_{n = 1}^{\I} Y_n = \{ y \},$
then $A_{ \{ y \} } = \dirlim A_{Y_n}.$
\end{rmk}

We need the following description of $C^* (\Z, X, h)_Y$ when
$h$ is minimal and $\sint (Y) \neq \varnothing.$
These results are originally from the unpublished preprint~\cite{Lnq}.
Outlines of the proofs can be found in
Section~3 of the survey article~\cite{LP1};
details (of this and much more) will appear in~\cite{LP3}.

The first theorem gives the decomposition of $X$ into Rokhlin
towers associated to $Y$ and the first return times to $Y.$
The bases of the towers are taken to be subsets of $h (Y)$ rather
than, as would be more usual, subsets of $Y.$
This is the convenient choice
for use with our definition of $C^* (\Z, X, h)_Y.$
For \rsha s and \rshd s as in the second theorem,
see Section~1 of~\cite{Ph6}.

\begin{thm}\label{Rokhlin}
Let $X$ be an infinite compact metric space,
and let $h \colon X \to X$ be a \mh.
Let $Y \subset X$ be closed and have nonempty interior.
For $y \in Y$ set $r (y) = \min \{ n \geq 1 \colon h^n (y) \in Y \}.$
Then $\sup_{y \in Y} r (y) < \infty.$
Let $n (0) < n (1) < \cdots < n (l)$ be the distinct
values in the range of $r,$ and for $0 \leq k \leq l$ set
\[
Y_k = {\overline{\{ y \in Y \colon r (y) = n (k) \} }}
\andeqn
Y_k^{\circ} = \sint (\{ y \in Y \colon r (y) = n (k) \}).
\]
Then:
\begin{itemize}
\item[(1)]
The sets $h^j (Y_k^{\circ}),$
for $0 \leq k \leq l$ and $1 \leq j \leq n (k),$ are disjoint.
\item[(2)]
$\bigcup_{k = 0}^l \bigcup_{j = 1}^{n (k)} h^j (Y_k) = X.$
\item[(3)]
$\bigcup_{k = 0}^l h^{n (k)} (Y_k) = Y.$
\end{itemize}
\end{thm}

\begin{thm}\label{SubalgThm}
Let $X,$ $h,$ $Y,$ and the other notation
be as in Notation~\ref{SubalgNtn} and Theorem~\ref{Rokhlin},
and let $A = C^* (\Z, X, h).$
Set
$B_Y = \bigoplus_{k = 0}^l  C (Y_k, M_{n (k)}).$
Then $A_Y$ has a \rshd\  with base spaces
$Y_0, Y_1, \ldots, Y_l$
and with standard representation $\sm \colon A_Y \to B_Y,$
such that:
\begin{itemize}
\item[(1)]
If $f \in C (X) \subset A_Y$ then
$\sm (f)_k \in C (Y_k, M_{n (k)})$ is given by
\[
\sm (f)_k = \diag ( (f |_{h (Y_k)}) \circ h, \,
                      (f |_{h^2 (Y_k)}) \circ h^2, \, \ldots, \,
                      (f |_{h^{n (k)} (Y_k)}) \circ h^{n (k)} ).
\]
\item[(2)]
If $f \in C_0 (X \setminus Y)$ then
$\sm (u f)_k \in C (Y_k, M_{n (k)})$ is given by
\[
\sm_k (u f)
 = \left( \begin{array}{ccccc}
  0      & 0      & \cdots & 0      & 0      \\
  1      & 0      & \cdots & 0      & 0      \\
  0      & 1      & \cdots & 0      & 0      \\
  \vdots & \vdots & \ddots & \vdots & \vdots \\
  0      & 0      & \cdots & 1      & 0      \end{array} \right)
 \sm_k (f).
\]
\end{itemize}
\end{thm}

The following result is due to Qing Lin
(Proposition~12 of~\cite{Lnq}).
Since that paper has not been published,
we give a proof here.

\begin{prp}\label{AySimple}
Let $X$ be an infinite compact metric space,
and let $h \colon X \to X$ be a \mh.
Let $A = C^* (\Z, X, h),$ and adopt Notation~\ref{SubalgNtn}.
Let $y \in X.$
Then $A_{ \{ y \} }$ is simple.
\end{prp}

\begin{proof}
Let $I \S A_{ \{ y \} }$ be a nonzero ideal.
Then $I \cap C (X)$ is an ideal in $C (X),$
so we can write $I \cap C (X) = C_0 (U)$ for some open set $U,$
which is necessarily given by
\[
U = \{ x \in X \colon {\mbox{there is $f \in I \cap C (X)$
            such that $f (x) \neq 0$}} \}.
\]

We first claim that $U \neq \E.$
Write $A_{ \{ y \} } = \dirlim A_{Y_m}$ as in Remark~\ref{LimSubalg},
with $Y_m$ chosen so that $\sint (Y_m) \neq \E.$
Then there exists $m$ such that $A_{Y_m} \cap I \neq \{ 0 \}.$
Let $a$ be a nonzero element of this intersection.
Using Theorem~\ref{SubalgThm}, one can fairly easily prove that
there is $N$ such that every element of $A_{Y_m}$ can be written
in the form $\sum_{n = - N}^N f_n u^n$ with $f_n \in C (X)$
for $- N \leq n \leq N.$
Moreover, if $a \neq 0$ and one writes
$a^* a = \sum_{n = - N}^N f_n u^n,$
then $f_0 \neq 0.$
Choose $x \in X$ such that $f_0 (x) \neq 0,$
choose a neighborhood $V$ of $x$ such that the sets $h^n (V),$
for $- N \leq n \leq N,$ are disjoint,
and choose $g \in C (X)$ such that $\supp (g) \S V$ and
$g (x) \neq 0.$
Then $g \in A_{ \{ y \} },$ and one checks that
\[
g a^* a g = \sum_{n = - N}^N g f_n u^n g
   = \sum_{n = - N}^N g (g \circ h^n) f_n u^n = g^2 f.
\]
So $g^2 f$ is a nonzero element of $I \cap C (X),$
proving the claim.

We next claim that
$h^{-1} ( U \SM \{ h (y) \}) \S U.$
So let $x \in U \SM \{ h (y) \}.$
Choose $f \in I \cap C (X)$ such that $f (x) \neq 0,$
and choose $g \in C_0 (X \SM \{ y \})$ such that $g (h^{-1} (x)) \neq 0.$
Then $u g \in A_{ \{ y \} },$
and
\[
(u g)^* f (u g) = {\overline{g}} u^* f u g = | g |^2 (f \circ h).
\]
Thus $| g |^2 (f \circ h) \in I \cap C (X)$ and is nonzero
at $h^{-1} (x).$
This proves the claim.

We further claim that $h ( U \SM \{ y \}) \S U.$
The proof is similar: let $x \in U \SM \{ y \},$
let $f \in I \cap C (X)$ and $g \in C_0 (X \SM \{ y \})$
be nonzero at $x,$ and consider $(u g) f (u g)^*.$

Now set $Z = X \SM U.$
The last two claims above imply that if $x \in X$ and $x$ is not in
the orbit of $y,$
then $h^k (x) \in Z$ for all $k \in \Z.$
Since $h$ is minimal, $Z$ is closed, and $Z \neq X,$
this is impossible.
If $h^n (y) \in Z$ for some $n > 0,$
then $h^{-1} ( U \SM \{ h (y) \}) \S U \SM \{ h (y) \}$
implies $h^k (y) \in Z$ for all $k \geq n.$
Since $\{ h^k (y) \colon k \geq n \}$ is dense by minimality,
this is also a contradiction.
Similarly, if $h^n (y) \in Z$ for some $n \leq 0,$
then $Z$ would contain the dense set
$\{ h^k (y) \colon k \leq n \},$ again a contradiction.
So $Z = \E,$ whence $U = X,$ and $1 \in I.$
\end{proof}

For convenience, we also reproduce the following
result of Qing Lin (Proposition~16 of~\cite{Lnq}).
The proof is sketched in the proof of Theorem~1.2 of~\cite{LP1}.

\begin{prp}\label{AyResTr}
Let $X$ be an infinite compact metric space,
and let $h \colon X \to X$ be a \mh.
Let $A = C^* (\Z, X, h),$ and adopt Notation~\ref{SubalgNtn}.
Let $y \in X.$
Then the restriction map $T (A) \to T (A_{ \{ y \} })$ is bijective.
\end{prp}

\section{The tracial rank of $A_{ \{ y \} }$}\label{Sec:TRSA}

Let the notation be as in the previous section,
with $A = C^* (\Z, X, h).$
We prove in this section that if the space $X$ is finite dimensional,
$h$ is minimal, and the image of $K_0 (A)$ in $\Aff (T (A))$ is dense,
then $A_{ \{ y \} }$ has tracial rank zero for ``most'' $y \in X.$
The result is surely true for all $y \in X,$
and the proof should not be much harder,
but this version is more convenient to prove
and suffices for our purposes.

If $A$ has at most countably many extreme tracial states,
equivalently, if $h$ has
at most countably many ergodic probability measures,
then this result (for all $y$) is easily
obtained from Theorem~4.4 of~\cite{Ph8}.
The point of the argument here is that it
gives the result when $T (A)$ is completely arbitrary.

\begin{dfn}\label{UNullD}
Let $X$ be a compact metric space,
and let $h \colon X \to X$ be a homeomorphism.
A Borel set $Z \subset X$ is called {\emph{universally null}}
if $\mu (Z) = 0$ for every $h$-invariant Borel probability measure
$\mu$ on $X.$
\end{dfn}

\begin{prp}\label{UNull}
Let $X$ be an infinite compact metric space
with finite covering dimension,
and let $h \colon X \to X$ be a \hme\  whose
set of periodic points has dimension zero or is empty.
For every $\ep > 0$ there are disjoint open sets
$U_1, U_2, \ldots, U_m \subset X$
such that $\diam (U_j) < \ep$ and $\sint ( {\overline{U_j}} ) = U_j$
for $1 \leq j \leq m,$
and such that $X \setminus \bigcup_{j = 1}^m U_j$ is universally null.
\end{prp}

\begin{proof}
Let $d = \dim (X).$
We apply Lemma~3.7 of~\cite{Kl}.
For the terminology used there, see Section~1 of the paper,
and note that $\om = \N.$
Taking $i$ there to be greater than $1 / \ep,$
we obtain closed subsets $T_1, T_2, \ldots, T_n \subset X$
which cover $X,$
whose interiors $U_1, U_2, \ldots, U_m$ are disjoint
and satisfy ${\overline{U_j}} = T_j,$
and such that whenever
\[
(j (0), \, n (0)), \, (j (1), \, n (1)), \, \ldots, \, (j (d), \, n (d))
  \in \{ 1, 2, \ldots, m \} \times \Z
\]
are $d + 1$ distinct pairs, then
\[
h^{n (0)} ( \partial T_{j (0)} ) \cap
h^{n (1)} ( \partial T_{j (1)} ) \cap \cdots \cap
h^{n (d)} ( \partial T_{j (d)} )
 = \varnothing.
\]

Set
\[
Z = X \setminus \bigcup_{j = 1}^m U_j
        \subset \bigcup_{j = 1}^m \partial T_j.
\]
Then $f = \sum_{n \in \Z} \ch_{h^n (Z)}$ is a nonnegative Borel function
on $X$ which is bounded by $d.$
If $\mu$ is any $h$-invariant Borel probability measure on $X,$
then
\[
\sum_{n \in \Z} \mu (h^n (Z)) = \int_X f \, d \mu \leq d.
\]
Since $\mu (h^n (Z)) = \mu (Z)$ for all $n,$
it follows that $\mu (Z) = 0.$
\end{proof}

\begin{cor}\label{UNullBd}
Let $X$ be an infinite compact metric space
with finite covering dimension,
and let $h \colon X \to X$ be a \mh.
Then there exists a dense $G_{\dt}$-set $S \subset X$ such that
for every $y \in S$ and every open set $V$ containing $X,$
there is a closed set $Y$ such that $y \in \sint (Y) \subset Y \subset V$
and $\partial Y$ is universally null.
\end{cor}

\begin{proof}
Apply Proposition~\ref{UNull} with $\ep = \frac{1}{n},$
and let $W_n$ be the union of the resulting finite collection
of open sets $U_{n, 1}, U_{n, 2}, \ldots, U_{n, m (n)}.$
Then $W_n$ is dense, since by minimality every
$h$-invariant Borel probability measure is nonzero on every nonempty
open subset of $X.$
So $S = \bigcap_{n = 1}^{\infty} W_n$ is a dense $G_{\dt}$-set in $X.$
Given $y \in S$ and an open set $V$ containing $y,$
for $n$ large enough we take $Y = {\overline{U_{n, j} }}$ with
$j$ chosen such that $y \in U_{n, j}.$
Note that $\partial Y \subset X \setminus W_n.$
\end{proof}

\begin{cor}\label{PartNullBd}
Let $X$ be an infinite compact metric space
with finite covering dimension,
and let $h \colon X \to X$ be a \mh.
Let $Y \subset X$ be a closed subset with nonempty interior
such that $\partial Y$ is universally null,
and adopt the notation of Theorem~\ref{Rokhlin}.
Then for $0 \leq k \leq l$ there are disjoint open subsets
$U_{k, 1}, U_{k, 2}, \dots, U_{k, r (k)}$ of $X$ which are
contained in $Y_l^{\circ},$
such that
\[
X \setminus \bigcup_{k = 0}^l \bigcup_{j = 1}^{n (k)}
             \bigcup_{i = 1}^{r (k)} h^j (U_{k, i})
\]
is universally null,
such that $\diam ( h^j (U_{k, i}) ) < \ep$
for $0 \leq k \leq l,$ $1 \leq j \leq n (k),$ and $1 \leq i \leq r (k),$
and such that the sets $h^j (U_{k, i}),$
for these $j, k, i,$ are disjoint.
\end{cor}

\begin{proof}
Choose $\dt > 0$ so small that if $d (x, y) < \dt$ then
\[
d (x, y) < \ep, \,\,\,\,\,\, d (h (x), h (y)) < \ep, \,\,\,\,\,\,
\ldots, \,\,\,\,\,\, d (h^{n (l)} (x), \, h^{n (l)} (y)) < \ep.
\]
Apply Proposition~\ref{UNull} with $\dt$ in place of $\ep.$
Then take $U_{k, 1}, U_{k, 2}, \dots, U_{k, r (k)}$ to be the
nonempty sets among
$U_1 \cap Y_k^{\circ}, \, U_2 \cap Y_k^{\circ},
               \, \ldots, \, U_m \cap Y_k^{\circ}.$
Note that
$\partial (U_i \cap Y_k^{\circ})
               \S \partial U_i \cup \partial Y_k^{\circ}.$
Disjointness follows from Theorem~\ref{Rokhlin}(1).
\end{proof}

The following lemma is taken from the proof of Lemma~3 of~\cite{Lnnaf}.
(It has also appeared in other places.)

\begin{lem}\label{NbhdNulln}
Let $A$ be a simple unital \ca, and let $B \subset A$ be a
separable commutative unital \csa,
which we identify with $C (X)$ for some compact metric space $X.$
For any $\ta \in T (A)$ let
$\mu_{\ta}$ be the Borel
probability measure on $X$ induced  by $\ta.$
Let $Z \subset X$ be a closed subset,
and suppose that $\mu_{\ta} (Z) = 0$ for all $\ta \in T(A).$
Then for any $\ep > 0$ there is an open subset $U \subset X$ containing
$Z$ such that
$\mu_{\ta} (U) < \ep$ for all $\ta \in T(A).$
\end{lem}

\begin{proof}
Choose open sets
$U_1 \supset {\overline{U}}_2 \supset U_2 \supset \cdots \supset Z$
such that $\bigcap_{n = 0}^{\infty} U_n = Z.$
Choose \cfn s $f_n \colon X \to [0, 1]$
such that $\supp (f_n) \subset U_n$ and $f_n = 1$ on $U_{n + 1}.$
Then $f_1 \geq f_2 \geq \cdots.$
Each $f_n$ defines a function $g_n \in \Aff (T (A))$
by $g_n (\ta) = \ta (f_n).$
Since $\mu_{\ta} (Z) = 0$ for all $\ta \in T(A),$ we have
$g_n (\ta) \to 0$ for each $\ta \in T (A)$
by the Dominated Convergence Theorem.
Since $0$ is continuous on the compact set $T (A),$
Dini's Theorem implies that $g_n \to 0$ uniformly on $T (A).$
Thus, for any $\ep > 0$ there exists $n$ such that
$\ta (f_n) < \ep$ for all $\ta \in T (A).$
It follows that $\mu_{\ta} (U_n) < \ep$ for all $\ta \in A.$
\end{proof}

\begin{cor}\label{NbhdNull}
Let $X$ be a compact metric space,
and let $h \colon X \to X$ be a \mh.
Let $Z \subset X$ be a universally null closed set,
and let $\ep > 0.$
Then there exists an open set $U \S X$ containing $Z$
such that $\mu (U) < \ep$
for every $h$-invariant Borel probability measure $\mu$ on $X.$
\end{cor}

\begin{proof}
Apply Lemma~\ref{NbhdNulln} with $A = C^*(\Z, X, h)$
and with $B$ being the canonically embedded copy of $C (X).$
The result follows from the one to one correspondence
between $h$-invariant Borel probability measures on $X$
and tracial states on $C^*(\Z, X, h).$
\end{proof}

\begin{thm}\label{TRZero}
Let $X$ be an infinite compact metric space
with finite covering dimension,
and let $h \colon X \to X$ be a \mh.
Let $A = C^* (\Z, X, h),$ and adopt Notation~\ref{SubalgNtn}.
Suppose that the image of
$K_0 (A)$ in $\Aff (T (A))$ is dense.
Then there is a dense $G_{\dt}$-set $S \subset X$ such that
for every $y \in S,$
the subalgebra $A_{ \{ y \} }$ has tracial rank zero.
\end{thm}

\begin{proof}
We begin with the following observations
(justified below), which hold for every $y \in X$:
\begin{itemize}
\item[(1)]
$A_{ \{ y \} }$ is simple.
\item[(2)]
$A_{ \{ y \} }$ has stable rank one.
\item[(3)]
The order on \pj s in
$A_{ \{ y \} }$ is determined by traces,
that is, if $p, q \in A_{ \{ y \} }$
are \pj s such that $\ta (p) < \ta (q)$ for every tracial state
$\ta,$ then $p \precsim q.$
\item[(4)]
The image of $K_0 (A_{ \{ y \} })$
in $\Aff (T (A_{ \{ y \} }))$ is dense.
\end{itemize}
Statement~(1) is Proposition~~\ref{AySimple}.
Statement~(2) is Theorem~4.2(1) of~\cite{Ph7},
and~(3) is Theorem~4.2(3) of~\cite{Ph7}.
For~(4), we observe that
the restriction map $T (A) \to T (A_{ \{ y \} })$
is bijective
(Proposition~\ref{AyResTr}),
and that $K_0 (A_{ \{ y \} }) \to K_0 (A)$
is a group isomorphism
(Theorem~4.1(3) of~\cite{Ph7}).
(The map $K_0 (A_{ \{ y \} }) \to K_0 (A)$ is an order isomorphism,
but we do not need this fact.)
Then~(4) follows from the hypothesis that the image of
$K_0 (A)$ in $\Aff (T (A))$ is dense.

Now let $S \S X$ be as in Corollary~\ref{UNullBd}, and let $y \in S.$
Given the properties above, we will prove that $A_{ \{ y \} }$ has tracial
rank zero by verifying the following condition:
For every finite subset ${\mathcal{F}} \S A$ and every $\ep > 0,$
there exists a \pj\  $p \in A$ and a unital finite dimensional
subalgebra $E \S p A p$ such that:
\begin{itemize}
\item[(5)]
$\| [a, p] \| < \ep$ for all $a \in {\mathcal{F}}.$
\item[(6)]
For all $a \in {\mathcal{F}}$ there exists $b \in E$
such that $\| p a p - b \| < \ep.$
\item[(7)]
$\ta (1 - p) < \ep$ for all tracial states $\ta$ on $A.$
\end{itemize}
That this suffices is Corollary~6.15 of~\cite{Ln2};
the Fundamental Comparison Property there
(stated in Theorem~6.8 of~\cite{Ln2}) is just~(3) above.

In fact, we need only let the finite set ${\mathcal{F}}$ run through
a cofinal collection of finite subsets of a set ${\mathcal{G}}$
which generates $A$ as a \ca.
Set
\[
{\mathcal{G}}_0 = \{ f \in C (X)
      \colon {\mbox{$f$ vanishes on a \nbhd\  of $y$}} \}.
\]
Then take
${\mathcal{G}} = \{ 1 \} \cup {\mathcal{G}}_0 \cup u {\mathcal{G}}_0.$
Our finite sets will have the form
${\mathcal{F}} = \{ 1 \} \cup {\mathcal{F}}_0 \cup u {\mathcal{F}}_0$
for a finite set ${\mathcal{F}}_0 \S {\mathcal{G}}_0.$

Use the hypothesis $y \in S$
to choose closed sets $Z_1 \supset Z_2 \supset \cdots$
such that $\bigcap_{m = 1}^{\I} Z_m = \{ y \},$
such that $y \in \sint (Z_m)$ for all $m,$
and such that $\partial Z_m$ is universally null for all $m.$
Then $A_{ \{ y \} } = \dirlim A_{Z_m}$ by Remark~\ref{LimSubalg}.
Given ${\mathcal{F}}_0 \S {\mathcal{G}}_0$ finite,
choose $m_0$ so large that $\supp (f) \S X \SM Z_{m_0}$
for all $f \in {\mathcal{F}}_0.$
Set $Y = Z_{m_0},$
and let the notation be as in Theorems~\ref{Rokhlin} and~\ref{SubalgThm}.
Choose $\dt > 0$ so small that whenever
$x_1, x_2 \in X$ satisfy $d (x_1, x_2) < \dt$ then
$| f (x_1) - f (x_2) | < \tfrac{1}{3} \ep$ for all $f \in {\mathcal{F}}.$
Use Corollary~\ref{PartNullBd} to find, for $0 \leq k \leq l,$
disjoint open subsets
$U_{k, 1}, U_{k, 2}, \dots, U_{k, r (k)}$ of $X$ which are
contained in $Y_l^{\circ},$
such that,
with
\[
U = \bigcup_{k = 0}^l \bigcup_{j = 1}^{n (k)}
             \bigcup_{i = 1}^{r (k)} h^j (U_{k, i}),
\]
the set $X \setminus U$ is universally null,
and such that $\diam ( h^j (U_{k, i}) ) < \dt$
for $0 \leq k \leq l,$ $1 \leq j \leq n (k),$ and $1 \leq i \leq r (k).$
By Corollary~\ref{NbhdNull}, there is an open set $V_0 \S X$ such that
$X \SM U \S V_0$ and $\mu (V_0) < \ep$
for every $h$-invariant Borel probability measure $\mu$ on $X.$
Choose an open set $V$ with
$X \SM U \S V \S {\overline{V}} \S V_0.$
Set $W_{k, i} = (X \SM {\overline{V}}) \cap U_{k, i}$
for $0 \leq k \leq l$ and $1 \leq i \leq r (k).$
Choose \cfn s
\[
g^{(0)}_{k, i}, \, g^{(1)}_{k, i}, \, g^{(2)}_{k, i}
  \colon X \to [0, 1]
\]
such that
\[
\supp ( g^{(2)}_{k, i} ) \S h (U_{k, i}), \,\,\,\,\,\,
g^{(2)}_{k, i} g^{(1)}_{k, i} = g^{(1)}_{k, i},
\]
\[
g^{(1)}_{k, i} g^{(0)}_{k, i} = g^{(0)}_{k, i}, \andeqn
g^{(0)}_{k, i} \ch_{h (W_{k, i})} = \ch_{h (W_{k, i})}
\]
for $0 \leq k \leq l$ and $1 \leq i \leq r (k).$

Because $A_{ \{ y \} } = \dirlim A_{Z_m}$ and using~(4),
Proposition~3.5 of~\cite{Ph8} provides \pj s
$q_{k, i} \in A_{ \{ y \} }$
(actually, in $A_{Z_m}$ for some $m \geq m_0$)
such that
\[
g^{(2)}_{k, i} q_{k, i} = q_{k, i}  \andeqn
q_{k, i} g^{(0)}_{k, i} = g^{(0)}_{k, i}
\]
for $0 \leq k \leq l$ and $1 \leq i \leq r (k).$
Then define $p_{k, i, j} = u^{j - 1} q_{k, i} u^{- (j - 1)}$
(so that $p_{k, i, 1} = q_{k, i}$),
and set
\[
p_{k, i} = \sum_{j = 1}^{n (k)} p_{k, i, j}
\andeqn
E_{k, i}
 = {\mathrm{span}}
       \left( \left\{ u^{j_1 - 1} q_{k, i} u^{j_2 - 1} \colon
                       1 \leq j_1, j_2 \leq n (k) \right\} \right).
\]

It is now convenient to introduce the index sets
\[
Q = \{ (k, i, j) \colon
  {\mbox{$0 \leq k \leq l,$ $1 \leq i \leq r (k),$ and
                $1 \leq j \leq n (k)$}} \}
\]
and
\[
Q_0 = \{ (k, i, j) \colon
  {\mbox{$0 \leq k \leq l,$ $1 \leq i \leq r (k),$ and
                $1 \leq j \leq n (k) - 1$}} \}.
\]
Thus, using Theorem~\ref{Rokhlin}(3) for the second equation,
\[
U = \bigcup_{(k, i, j) \in Q} h^j (U_{k, i})
\andeqn
U \cap (X \SM Y) = \bigcup_{(k, i, j) \in Q_0} h^j (U_{k, i}).
\]

We claim that (proofs below):
\begin{itemize}
\item[(8)]
$p_{k, i} \in A_{ \{ y \} }$ and $E_{k, i} \S A_{ \{ y \} }.$
\item[(9)]
$(g^{(2)}_{k, i} \circ h^{- (j - 1)}) p_{k, i, j} = p_{k, i, j}$
and
$p_{k, i, j} (g^{(0)}_{k, i} \circ h^{- (j - 1)})
  = g^{(0)}_{k, i} \circ h^{- (j - 1)}$
for $(k, i, j) \in Q.$
\item[(10)]
The $p_{k, i, j},$ for $(k, i, j) \in Q,$ are orthogonal \pj s.
\item[(11)]
The $p_{k, i},$ for $0 \leq k \leq l$ and $1 \leq i \leq r (k),$
are orthogonal \pj s.
\item[(12)]
$E_{k, i}$ is a finite dimensional \ca\  with identity $p_{k, i}$
(which is isomorphic to $M_{n (k)},$
and has matrix units $e_{j_1, j_2} = u^{j_1 - 1} q_{k, i} u^{j_2 - 1}$).
\end{itemize}
Then we set
\[
p = \sum_{k = 0}^{l} \sum_{i = 1}^{r (k)} p_{k, i}
\andeqn
E = \bigoplus_{k = 0}^{l} \bigoplus_{i = 1}^{r (k)} p_{k, i}.
\]
It follows from the claims above that:
\begin{itemize}
\item[(13)]
$p$ is a \pj\  in $A_{ \{ y \} }.$
\item[(14)]
$E$ is a finite dimensional \csa\  of $A_{ \{ y \} }$ with identity $p.$
\end{itemize}
We further claim that (again, proofs below):
\begin{itemize}
\item[(15)]
$\| [a, p] \| < \ep$ for all $a \in {\mathcal{F}}.$
\item[(16)]
For all $a \in {\mathcal{F}}$ there exists $b \in E$
such that $\| p a p - b \| < \ep.$
\item[(17)]
$\ta (1 - p) < \ep$ for all tracial states $\ta$ on $A_{ \{ y \} }.$
\end{itemize}

Statements~(13) through~(17) yield the condition for tracial rank zero
given above, so the proof will be complete once~(8) through~(12)
and~(15) through~(17) are verified.

We begin with~(8).
It is enough to prove that $u^{j - 1} q_{k, i} \in A_{ \{ y \} }$
for $1 \leq j \leq n (k).$
We write
\[
u^{j - 1} q_{k, i}
  = u^{j - 1} (g^{(2)}_{k, i})^{j - 1} q_{k, i}
  = u ( g^{(2)}_{k, i} \circ h^{- (j - 2)}) \cdots
        u ( g^{(2)}_{k, i} \circ h^{- 1}) \cdot
        u g^{(2)}_{k, i} \cdot q_{k, i}.
\]
The functions
$g^{(2)}_{k, i}, \, g^{(2)}_{k, i} \circ h^{- 1}, \, \ldots,
 \, g^{(2)}_{k, i} \circ h^{- (j - 2)}$
are supported in
\[
h (U_{k, i}) \S h (Y_k), \, h^2 (U_{k, i}) \S h^2 (Y_k), \, \ldots,
   h^{j - 1} (U_{k, i}) \S h^{j - 1} (Y_k),
\]
so that every factor except $q_{k, i}$ in the
last expression is in $u C_0 (X \setminus Y).$
This proves~(8).

The relations in~(9) are obtained by conjugating the relations
\[
g^{(2)}_{k, i} q_{k, i} = q_{k, i}  \andeqn
q_{k, i} g^{(0)}_{k, i} = g^{(0)}_{k, i}
\]
by $u^{j - 1}$ and using
$u^{(j - 1)} g u^{- (j - 1)} = g \circ h^{- (j - 1)}.$

Part~(10) follows from~(9) together with the observation that
$\supp (g^{(2)}_{k, i} \circ h^{- (j - 1)}) \S h^j (U_{k, i})$
and the fact that the sets $h^j (U_{k, i}),$ for $(k, i, j) \in Q,$
are disjoint.
Part~(11) is now immediate.
In~(12), the matrix unit relations are easy once
one has the orthogonality relations in~(10),
and that
$E_{k, i}$ is a finite dimensional \ca\  with identity $p_{k, i}$
is then immediate.

We now prove~(15) and~(16).
By construction, for $(k, i, j) \in Q$
there exists $\ld_{k, i, j} \in \C$ such that
$| \ld_{k, i, j} - f (x) | < \tfrac{1}{3} \ep$
for every $f \in {\mathcal{F}}_0$ and every $x \in h^j (U_{k, i}).$
Since $f$ vanishes on $X \SM Y,$
we may require $\ld_{k, i, j} = 0$ when
$(k, i, j) \not\in Q_0,$ that is, when $j = n (k).$
Set
\[
b = \sum_{(k, i, j) \in Q} \ld_{k, i, j} p_{k, i, j},
\]
which is in $E.$

We have
\[
\left\| f \cdot \ts{\ds{\sum}_{(k, i, j) \in Q}}
                  g^{(2)}_{k, i} \circ h^{- (j - 1)}
    - \ts{\ds{\sum}_{(k, i, j) \in Q}}
             \ld_{k, i, j} \cdot (g^{(2)}_{k, i} \circ h^{- (j - 1)}) \right\|
 \leq \frac{\ep}{3}.
\]
Multiplying on the right by
\[
p = \sum_{(k, i, j) \in Q} p_{k, i, j},
\]
using~(9) and disjointness of the supports of the
$g^{(2)}_{k, i} \circ h^{- (j - 1)}$ for $(k, i, j) \in Q,$
we get $\left\| f p - b \right\| \leq \tfrac{1}{3} \ep.$
Similarly
$\left\| p f - b \right\| \leq \tfrac{1}{3} \ep$
and $\left\| p f p - b \right\| \leq \tfrac{1}{3} \ep < \ep.$
Thus $\| [p, f] \| \leq \tfrac{2}{3} \ep < \ep.$
This gives~(15) and~(16) for $a \in {\mathcal{F}}_0.$

Next, set
\[
e = \sum_{(k, i, j) \in Q_0} p_{k, i, j} \leq p.
\]
Because $\ld_{k, i, j} = 0$ when $j = n (k),$
we get $e b = b e = b,$ and also
\[
u b = \sum_{(k, i, j) \in Q_0}
             \ld_{k, i, j} \cdot u p_{k, i, j}
    = \sum_{(k, i, j) \in Q_0}
             \ld_{k, i, j} \cdot u^{j} q_{k, i} u^{- (j - 1)} \in E.
\]
The next step is to show that $p u f = u e f.$
Now
\[
u^* p u - e = \sum_{k = 0}^{l} \sum_{i = 1}^{r (k)} u^* q_{k, i} u.
\]
We have
\[
(u^* q_{k, i} u) (g^{(2)}_{k, i} \circ h)
  = u^* q_{k, i} g^{(2)}_{k, i} u = u^* q_{k, i} u,
\]
and because $\supp (g^{(2)}_{k, i} \circ h) \S Y$
and $f$ vanishes on $Y,$ we get $(g^{(2)}_{k, i} \circ h) f = 0.$
Therefore
\[
p u f - u e f = u (u^* p u - e) f
  = u \sum_{k = 0}^{l} \sum_{i = 1}^{r (k)} u^* q_{k, i} u
                 (g^{(2)}_{k, i} \circ h) f
  = 0.
\]
Now
\[
\| p u f - u b \| = \| e f - b \|
   \leq \| e \| \cdot \| p f - b \|
   \leq \tfrac{1}{3} \ep
\andeqn
\| u f p - u b \| = \| f p - b \| \leq \tfrac{1}{3} \ep,
\]
so
\[
\| [p, u f] \| \leq \tfrac{2}{3} \ep < \ep
\andeqn
\| p u f p - u b \|
 \leq \| p u f - u b \| \cdot \| p \|
 \leq \tfrac{1}{3} \ep
 < \ep.
\]
This gives~(15) and~(16) for $a \in u {\mathcal{F}}_0.$
The only other element of ${\mathcal{F}}$ is $1,$ for which~(15) and~(16) are obvious.
This completes the verification of~(15) and~(16).

Finally, we prove~(17).
Let $\ta$ be a tracial state on $A_{ \{ y \} }.$
The restriction map $T (A) \to T (A_{ \{ y \} })$ is bijective
by Proposition~\ref{AyResTr}.
Therefore $\ta$ is the restriction to $A_{ \{ y \} }$
of the tracial state on $A$ obtained from some
$h$-invariant Borel probability measure $\mu.$
We have
$p \geq \sum_{(k, i, j) \in Q} g^{(0)}_{k, i} \circ h^{- (j - 1)}.$
Therefore
\begin{align*}
\ta (p)
  & \geq \sum_{(k, i, j) \in Q}
             \int_X g^{(0)}_{k, i} \circ h^{- (j - 1)} \, d \mu
    \geq \sum_{(k, i, j) \in Q} \mu ( h^j (W_{k, i}) )   \\
  & = \mu ((X \SM {\overline{V}}) \cap U)
    = \mu (X \SM {\overline{V}})
    \geq 1 - \mu (V_0)
    > 1 - \ep.
\end{align*}
This proves~(17),
and completes the proof of the theorem.
\end{proof}

\section{The tracial rank of the of the transformation group
                            C*-algebra}\label{Sec:TRCrP}

In this section, we prove the main result,
Theorem~\ref{TM} below.

The following lemma is a version of a result of L.~G.\  Brown.

\begin{lem}\label{L2}
Let $A$ be a \ca\  with real rank zero,
and let $a_1, a_2, a_3 \in A$ satisfy
\[
0 \leq a_1 \leq a_2 \leq a_3 \leq 1,
\,\,\,\,\,\,
a_3 a_2 = a_2, \andeqn a_2 a_1 = a_1.
\]
Then there exists a projection $e \in A$ such that
\[
e a_1 = a_1 e = a_1 \andeqn a_3 e = e a_3 = e.
\]
\end{lem}

\begin{proof}
Since the unitization of $A$ again has real rank zero,
and since the equation $a_3 e = e$ will imply $e \in A,$
we may assume $A$ is unital.
In $A^{**}$ let $p$ and $q$ be the closed \pj s
\[
p = \ch_{[2/3, \, \I]} (a_2)
\andeqn q = \ch_{[0, \, 1/3]} (a_2).
\]
Then one checks that
\[
a_3 (1 - q) = 1 - q \andeqn p a_1 = a_1.
\]
Use Theorem~1 of~\cite{Br}
to find a projection $e\in A$ such that $p \leq e \leq 1- q.$
This \pj\  clearly satisfies the conclusion.
\end{proof}

\begin{lem}\label{L3}
Let $X$ be an infinite compact metric space,
and let $h \colon X \to X$  be a minimal homeomorphism.
Let $A = C^* (\Z, X, h),$ and adopt Notation~\ref{SubalgNtn}.
Let $y \in X,$ and suppose that $A_{ \{ y \} }$ has
real rank zero and stable rank one.
Then for any
$\ep > 0$ and any finite subset ${\mathcal{F}} \subset A,$ there
is a projection $p \in  A_{ \{ y \} }$  such that:
\begin{itemize}
\item[(1)]
$\| p a - a p \| < \ep$ for all $a \in {\mathcal{F}}.$
\item[(2)]
$p a p \in p A_{ \{ y \} } p$ for all $a \in {\mathcal{F}}.$
\item[(3)]
$\ta (1 - p) < \ep$ for all $\ta \in T(A).$
\end{itemize}
\end{lem}

\begin{proof}
We may assume that ${\mathcal{F}} = {\mathcal{G}} \cup \{u\}$
for some finite subset ${\mathcal{G}} \subset C (X).$

Choose $N_0 \in \N$ so large that $4 \pi / N_0 < \ep.$
Choose $\dt_0 > 0$ with $\dt_0 < \frac{1}{2} \ep$ and so small
that $d (x_1, x_2) < 4 \dt_0$ implies
$| f (x_1) - f (x_2)| < \frac{1}{4} \ep$ for all $f \in {\mathcal{G}}.$
Choose $\dt > 0$ with $\dt \leq \dt_0$ and such that
whenever $d (x_1, x_2) < \dt$ and $0 \leq n \leq N_0,$ then
$d (h^{-n} (x_1), \, h^{-n} (x_2)) < \dt_0.$

Since $h$ is  minimal, there is $N > N_0 + 1$ such that
$d (h^N (y), \, y) < \dt.$
Choose $R \in \N$ so large that $R > (N + N_0 + 1) / \min (1, \ep).$
Since $h$ is free, there is an
open neighborhood $U$ of $y$ in $X$ such that
\[
h^{- N_0} (U), \, h^{- N_0 + 1} (U), \, \ldots,
\, U, \, h (U), \, \ldots, \, h^R (U)
\]
are disjoint and all have diameter less than $\dt.$
In particular, this is true with $N$ in place of $R,$
and furthermore $\mu (U) < \ep / (N + N_0 + 1)$
for every $h$-invariant Borel probability measure $\mu.$

Choose \cfn s $g_0, g_1, g_2, f_0 \colon X \to [0, 1]$ such that
\[
g_0 (y) = 1, \,\,\,\,\,\,
g_1 g_0 = g_0, \,\,\,\,\,\,
g_2 g_1 = g_1, \,\,\,\,\,\,
f_0 g_2 = g_2, \andeqn
\supp (f_0) \S U.
\]
By Lemma~\ref{L2},
there exists a \pj\  $q_0 \in A_{ \{ y \} }$ such that
\[
g_1 q_0 = q_0 g_1 = g_1 \andeqn f_0 q_0 = q_0 f_0 = q_0.
\]
For $- N_0 \leq n \leq N$ set
\[
q_n = u^n q_0 u^{- n}, \,\,\,\,\,\,
f_n = u^n f_0 u^{- n} = f_0 \circ h^{- n},
\andeqn U_n = h^n (U).
\]
We have $\supp (f_n) \S U_n$ and $f_n q_n = q_n f_n = q_n,$
so the $q_n$ are \mops\  in $A.$

We claim that $q_n \in A_{ \{ y \} }$ for $- N_0 \leq n \leq N.$
We first consider $q_{- n}$ for $1 \leq n \leq N_0.$
For these $n,$ we have $u f_{-n} \in A_{ \{ y \} },$ so
\begin{equation}\label{Eq:an_Dfn}
a_n = f_0^n u^n
    = (u f_{- 1}) (u f_{- 2}) \cdots (u f_{- n})
    \in A_{ \{ y \} }.
\end{equation}
Then $q_{-n} = u^{-n} q_0 u^n = a_n^* q_0 a_n \in A_{ \{ y \} },$
as desired.

Next, we show $q_1 \in A_{ \{ y \} }.$
{}From $q_0 g_0 = q_0 g_1 g_0 = g_1 g_0$ we get
$(q_0 - g_1) (f_0 - g_0) = q_0 - g_1.$
Also, $u (f_0 - g_0) \in A_{ \{ y \} }$ because $(f_0 - g_0) (y) = 0.$
Moreover, $u g_1 u^* = g_1 \circ h^{- 1} \in A_{ \{ y \} }.$
So
\[
q_1 = u q_0 u^*
    = u (q_0 - g_1) u^* + u g_1 u^*
    = [u (f_0 - g_0)] (q_0 - g_1) [u (f_0 - g_0)]^* + u g_1 u^*
    \in A_{ \{ y \} }.
\]
For $2 \leq n \leq N$ we now have
\[
q_n = [ (u f_n) (u f_{n - 1}) \cdots (u f_1) ] q_1
            [ (u f_n) (u f_{n - 1}) \cdots (u f_1) ]^*,
\]
which is in $A_{ \{ y \} }$ because $f_1, f_2, \ldots, f_N$
all vanish on $U$ and $y \in U.$
This proves the claim.

We now have a sequence of \pj s, in principle going to infinity
in both directions:
\[
\ldots, \, q_{- N_0}, \, \ldots, \, q_{-1}, \, q_0, \, q_1, \,
\ldots, \, q_{N - N_0}, \, \ldots, \, q_{N - 1}, \, q_N, \, \ldots.
\]
The ones shown are orthogonal, and conjugation by $u$ is the shift.
The \pj s $q_0$ and $q_N$ live over open sets which are disjoint
but close to each other, and similarly for the pairs
$q_{-1}$ and $q_{N - 1}$ down to $q_{- N_0}$ and $q_{N - N_0}.$
We are now going to use Berg's technique~\cite{Bg}
to splice this sequence along the pairs of indices
$(- N_0, \, N - N_0)$ through $(0, N),$
obtaining a loop of length $N$
on which conjugation by $u$ is approximately the cyclic shift.

We claim that there is a partial isometry $w \in  A_{ \{ y \} }$
such that $w^* w = q_0$ and $w w^* = q_N.$
Certainly $[q_0] = [q_N]$ in $K_0 (A).$
The map $K_0 (A_{ \{ y \} }) \to K_0 (A)$
is a group isomorphism by Theorem~4.1(3) of~\cite{Ph7},
so $[q_0] = [q_N]$ in $K_0 (A_{ \{ y \} }).$
Since $A_{ \{ y \} }$ has stable rank one,
Proposition~6.5.1 of~\cite{Bl} implies that \pj s in matrix algebras
over $A_{ \{ y \} }$ satisfy cancellation.
The claim follows.

For $t \in \R$ define
$v (t) = \cos (\pi t / 2) (q_0 + q_N) + \sin (\pi t / 2) (w - w^*).$
Then $v (t)$ is a unitary in the corner
$(q_0 + q_N) A_{ \{ y \} } (q_0 + q_N)$ whose matrix with respect to
the obvious block decomposition is
\[
v (t)
 = \left( \begin{array}{cc} \cos (\pi t / 2) & - \sin (\pi t / 2) \\
            \sin (\pi t / 2) & \cos (\pi t / 2) \end{array} \right).
\]
For $0 \leq k \leq N_0$ define $w_k = u^{- k} v (k / N_0) u^k.$
With $a_k$ as in~(\ref{Eq:an_Dfn}) and with
\[
b_k = f_N^k u^k
    = (u f_{N - 1}) (u f_{N - 2}) \cdots (u f_{N - k})
    \in A_{ \{ y \} }
\]
(because $N_0 < N$),
and using $U_{- k} \cap U_{N - k} = U_0 \cap U_{N} = \E,$
we can write
\[
w_k = (a_k + b_k)^* v (k / N_0) (a_k + b_k) \in A_{ \{ y \} }.
\]
Therefore $w_k$ is a unitary in the corner
$(q_{- k} + q_{N - k}) A_{ \{ y \} } (q_{- k} + q_{N - k}).$
Moreover,
adding estimates on the differences of the
matrix entries at the second step,
\[
\| u w_{k + 1} u^* - w_k \|
 = \| v (k / N_0) - v ((k - 1) / N_0) \|
 \leq 2 \pi / N_0
 < \tfrac{1}{2} \ep.
\]

Now define $e_n = q_n$ for $0 \leq n \leq N - N_0,$
and for $N - N_0 \leq n \leq N$ write $k = N - n$ and set
$e_n = w_k q_{- k} w_k^*.$
The two definitions for $n = N - N_0$ agree
because $w_{N_0} q_{- N_0} w_{N_0}^* = q_{N - N_0},$
and moreover $e_N = e_0.$
Therefore $u e_{n - 1} u^* = e_n$ for $1 \leq n \leq N - N_0,$
and also $u e_N u^* = e_1,$
while for $N - N_0 < n \leq N$ we have
\[
\| u e_{n - 1} u^* - e_n \|
 \leq 2 \| u w_{N - n + 1} u^* - w_{N - n} \|
 < \ep.
\]
Also, clearly $e_n \in A_{ \{ y \} }$ for all $n.$

Set $e = \sum_{n = 1}^{N} e_n$ and $p = 1 - e.$
We verify that $p$ satisfies~(1) through~(3).

First,
\[
p - u p u^*
 = u e u^* - e
 = \sum_{n = N_0 + 1}^{N} (u e_{n - 1} u^* - e_n ).
\]
The terms in the sum are orthogonal and have norm less than $\ep,$
so $\| u p u^* - p \| < \ep.$
Furthermore,
since
\[
1 - g_1 \in C_0 (X \SM \{ y \}), \,\,\,\,\,\, q_0 g_1 = g_1,
\andeqn p \leq 1 - q_0,
\]
we get
\[
p u p = p u (1 - g_1) (1 - q_0) p \in A_{ \{ y \} }.
\]
This is~(1) and~(2) for the element $u \in {\mathcal{F}}.$

Next, let $g \in {\mathcal{G}}.$
The sets $U_0, U_1, \ldots, U_N$ all have diameter less than $\dt.$
We have $d (h^N (y), \, y) < \dt,$
so the choice of $\dt$ implies that
$d (h^n (y), \, h^{n - N} (y)) < \dt_0$ for $N - N_0 \leq n \leq N.$
Also, $U_{n - N} = h^{n - N} (U_0)$
has diameter less than $\dt.$
Therefore $U_{n - N} \cup U_n$ has diameter less than
$2 \dt + \dt_0 \leq 3 \dt_0.$
Since $g$ varies by at most $\frac{1}{4} \ep$ on any set with
diameter less than $4 \dt_0,$
and since the sets
\[
U_1, \, U_2, \, \ldots, \, U_{N - N_0 - 1}, \,
U_{N - N_0} \cup U_{- N_0}, \,
U_{N - N_0 + 1} \cup U_{- N_0 + 1}, \, \ldots, \,
U_N \cup U_0
\]
are disjoint,
there is ${\widetilde{g}} \in C (X)$
which is constant on each of these sets and satisfies
$\| g - {\widetilde{g}} \| < \frac{1}{2} \ep.$
Let the values of ${\widetilde{g}}$ on these sets be
$\ld_1$ on $U_1$ through $\ld_N$ on $U_N \cup U_0.$

For $0 \leq n \leq N - N_0$ we have
$f_n e_n = f_n q_n = q_n = e_n,$ whence $\supp (f_n) \S U_n$ implies
\[
{\widetilde{g}} e_n
 = {\widetilde{g}} f_n e_n = \ld_n f_n e_n = \ld_n e_n
\]
and similarly $e_n {\widetilde{g}} = \ld_n e_n.$
For $N - N_0 < n \leq N$ we use
$e_n \in (q_{n - N} + q_{n}) A_{ \{ y \} } (q_{n - N} + q_{n})$
to get, in the same way as above,
\[
(f_{n - N} + f_{n}) e_n = e_n (f_{n - N} + f_{n}) = e_n,
\]
so from $\supp (f_{n - N} + f_{n}) \S U_{n - N} + U_{n}$ we get
${\widetilde{g}} e_n = \ld_n e_n = e_n {\widetilde{g}}.$
Since $\| g - {\widetilde{g}} \| < \frac{1}{2} \ep,$
it follows that
\[
\| p g - g p \| = \| g e - e g \| < \ep.
\]
This is~(1) for $g.$
That $p g p \in A_{ \{ y \} }$ follows from the fact that
$g$ and $p$ are in this subalgebra.
So we also have~(2) for~$g.$

It remains only to verify~(3).
Let $\ta \in T (A),$
and let $\mu$ be the corresponding $h$-invariant probability measure
on $X.$
We have
\[
1 - p = e \leq \sum_{n = - N_0}^N q_n \leq \sum_{n = - N_0}^N f_n,
\]
so
\[
\ta (1 - p) \leq \sum_{n = - N_0}^N \mu (U_n)
   = \sum_{n = - N_0}^N \mu (h^n (U))
   = (N + N_0 + 1) \mu (U)
   < \ep.
\]
This completes the proof.
\end{proof}

The next two lemmas assert that for a simple \ca,
in Definition~\ref{DTR0} of tracial rank zero,
and in the definition of
the local approximation property of Popa~\cite{Pp},
the finite dimensional subalgebras
can be replaced by subalgebras with tracial rank zero.
We state these lemmas separately so as to have them available
for use elsewhere.
We do not require the subalgebra $B$ in the first to be simple,
because it does not simplify the proof.
However, for the application in this paper,
$B$ will in fact be simple, so that the reader
need not worry about the extra complexity
of tracial rank for nonsimple algebras.

The first lemma is certainly known.
The second is a special case of Theorem~4.6 of~\cite{HLX},
where neither $A$ nor $B$ is required to be simple,
but the proof here is much easier.
This argument can also be found in the last part of the
proof of Theorem~2.9 of~\cite{LO}.

\begin{lem}\label{Popa}
Let $A$ be a simple unital \ca.
Suppose that for every $\ep > 0$
and every finite subset ${\mathcal{F}} \subset A,$
there exists a unital \csa\  $B \subset A$ which has tracial rank zero
and a projection $p \in B$ such that
\[
\| p a - a p \| < \ep \andeqn
{\mathrm{dist}} (p a p, \, p B p) < \ep
\]
for all $a \in {\mathcal{F}}.$
Then $A$ has the local approximation property of Popa~\cite{Pp},
that is,
for every $\ep > 0$
and every finite subset ${\mathcal{F}} \subset A,$
there exists a nonzero projection $q \in A$
and a finite dimensional unital \csa\  $D \subset q A q$ such that
\[
\| q a - a q \| < \ep \andeqn
{\mathrm{dist}} (q a q, \, D) < \ep
\]
for all $a \in {\mathcal{F}}.$
\end{lem}

\begin{proof}
Let $\ep > 0$ and let ${\mathcal{F}} \subset A$ be a finite subset.
Choose $B$ and $p$ as in the hypotheses,
with $\frac{1}{5} \ep$ in place of $\ep.$
In particular, there exists a finite subset
${\mathcal{G}} \subset p B p$ such that
${\mathrm{dist}} (p a p, \, {\mathcal{G}}) < \frac{1}{5} \ep$
for all $a \in {\mathcal{F}}.$

Since $B$ has tracial rank zero,
Theorem~5.3 of~\cite{Ln2}
implies that $p B p$ also has tracial rank zero.
So there is a nonzero \pj\  $q \in p B p$ and a
finite dimensional unital \csa\  $D \subset q B q$ such that
\[
\| q b - b q \| < \tfrac{1}{5} \ep
\andeqn
{\mathrm{dist}} (q b q, D) < \tfrac{1}{5} \ep
\]
all $b \in {\mathcal{G}}.$
We claim that
\[
\| q a - a q \| < \ep \andeqn
{\mathrm{dist}} (q a q, \, D) < \ep
\]
for all $a \in {\mathcal{F}}.$
The second part is easy.
For the first, let $a \in {\mathcal{F}}$
and choose $b \in {\mathcal{G}}$
such that $\| q a q - b \| < \tfrac{1}{5} \ep.$
Then
\begin{align*}
\| q a - a q \|
& = \| q p a - a p q \|
  \leq \| q p a p - p a p q \| + 2 \| p a - a p \|  \\
& \leq \| q b - b q \| + 2 \| b - p a p \| + 2 \| p a - a p \|
  < \ep.
\end{align*}
This completes the proof.
\end{proof}

\begin{lem}\label{TTR0}
Let $A$ be a simple unital \ca.
Suppose that
for every finite subset ${\mathcal{F}} \S A,$ every $\ep > 0,$
and every nonzero positive element $c \in A,$
there exists a \pj\  $p \in A$ and a simple unital
subalgebra $B \S p A p$ with tracial rank zero such that:
\begin{itemize}
\item[(1)]
$\| [a, p] \| < \ep$ for all $a \in {\mathcal{F}}.$
\item[(2)]
$\dist (p a p, B) < \ep$ for all $a \in {\mathcal{F}}.$
\item[(3)]
$p$ is Murray-von Neumann equivalent to a \pj\  in ${\overline{c A c}}.$
\end{itemize}
Then $A$ has tracial rank zero.
\end{lem}

\begin{proof}
Even without~(3), the algebra $A$ has Property~(SP) by
Lemma~2.12 of~\cite{Ln1}.
If $A$ is finite dimensional, the result is trivial;
otherwise,
there are nonzero orthogonal \pj s $e, f \in {\overline{c A c}}.$

The proof is now essentially the same as that of Lemma~\ref{Popa}.
The \pj\  $p$ in that proof can be required to satisfy
$1 - p \precsim e.$
Use Lemma~3.1 of~\cite{Ln1} to find a nonzero \pj\  $f_0 \leq p$
with $f_0 \precsim f,$ and require the \pj\  $q \in p B p$
to satisfy $p - q \precsim f_0$ in $p B p.$
\end{proof}

\begin{thm}\label{T1}
Let $X$ be an infinite compact metric space,
and let $h \colon X \to X$  be a minimal homeomorphism.
Let $A = C^* (\Z, X, h),$ and adopt Notation~\ref{SubalgNtn}.
Suppose that there is $y \in X$ such that $A_{ \{ y \} }$ has
tracial rank zero.
Then $A$ has tracial rank zero.
\end{thm}

\begin{proof}
We verify the condition of Lemma~\ref{TTR0}.
Thus, let ${\mathcal{F}} \S A$ be a finite subset,
let $\ep > 0,$
and let $c \in A$ be a nonzero positive element.
Beyond Lemma~\ref{L3},
the main step of the proof is to find a nonzero \pj\  in $A_{ \{ y \} }$
which is Murray-von Neumann equivalent
to a \pj\  in ${\overline{c A c}}.$

The algebra $A$ is simple,
so Lemma~\ref{L3}, Lemma~\ref{Popa}, and Lemma~2.12 of~\cite{Ln1}
imply that $A$ has property~(SP).
Therefore there is a nonzero projection $e \in {\overline{c A c}}.$
Set
\[
\dt_0 = \frac{1}{18} \inf_{\ta \in T (A)} \ta (e) \leq \frac{1}{18}.
\]
By Lemma~\ref{L3},
there is a projection $q \in  A_{ \{ y \} }$
and an element $b_0 \in q A_{ \{ y \} } q$ such that
\[
\| q e - e q \| < \dt_0, \,\,\,\,\,\,
\| q e q - b_0 \| < \dt_0,
\andeqn
\sup_{\ta \in T (A)} \ta (1 - q) \leq \dt_0.
\]
Replacing $b_0$ by $\frac{1}{2} (b_0 + b_0^*),$
we may assume that $b_0$ is selfadjoint.
We have $- \dt_0 \leq b_0 \leq 1 + \dt_0,$
so applying \ct\  functional calculus
we may find $b \in q A_{ \{ y \} } q$
such that $0 \leq b \leq 1$ and $\| q e q - b \| < 2 \dt_0.$
Using $\| q e - e q \| < \dt_0$ on the last term in the second
expression,
we get
\[
\| b^2 - b \| \leq 3 \| b - q e q \| + \| (q e q)^2 - q e q \|
     < 3 \cdot 2 \dt_0 + \dt_0 = 7 \dt_0 < \tfrac{1}{4}.
\]
Therefore there is a \pj\  $e_1 \in q A_{ \{ y \} } q$
such that $\| e_1 - b \| < 14 \dt_0,$
giving $\| e_1 -  q e q \| < 16 \dt_0.$
Similarly (actually, one gets a better estimate)
there is a \pj\  $e_2 \in (1 - q) A (1 - q)$
such that $\| e_2 - (1 - q) e (1 - q) \| < 16 \dt_0.$
Therefore
\[
\| e_1 + e_2 - [q e q + (1 - q) e (1 - q)] \| < 16 \dt_0
\]
and, using $\| q e - e q \| < \dt_0$ again, we have
$\| e_1 + e_2 - e \| < 18 \dt_0 \leq 1.$
It follows that $e_1 \precsim e.$
Also, for $\ta \in T (A),$ we have
\[
\ta (e_1)
 > \ta (q e q) - 16 \dt_0
 = \ta (e) - \ta ( (1 - q) e (1 - q) ) - 16 \dt_0
 \geq \ta (e) - \ta (1 - q) - 16 \dt_0
 > 0,
\]
so $e_1 \neq 0.$

Now set $\ep_0 = \inf ( \{\ta (e_2) \colon \ta \in T (A) \}).$
By Lemma~\ref{L3}, there is a
projection $p \in  A_{ \{ y \} }$  such that
\begin{itemize}
\item[(1)]
$\| p a - a p \| < \ep$ for all $a \in {\mathcal{F}}.$
\item[(2)]
$p a p \in p A_{ \{ y \} } p$ for all $a \in {\mathcal{F}}.$
\item[(3)]
$\ta (1 - p) < \ep_0$ for all $\ta \in T (A).$
\end{itemize}
Using Proposition~\ref{AyResTr},
it follows from~(3) that
$\ta (1 - p) < \ep_0$ for all $\ta \in T (A_{ \{ y \} }).$
Since $A_{ \{ y \} }$ has tracial rank zero,
and since tracial rank zero implies that the order on \pj s is
determined by traces (Theorems~6.8 and~6.13 of~\cite{Ln2}),
it follows that $1 - p \precsim e_1 \precsim e.$
Since $p A_{ \{ y \} } p$ also has
tracial rank zero (Theorem~5.3 of~\cite{Ln2}),
we have verified the hypotheses of Lemma~\ref{TTR0}.
Thus $A$ has tracial rank zero.
\end{proof}

\begin{thm}\label{TM}
Let $X$ be an infinite compact metric space
with finite covering dimension,
and let $h \colon X \to X$ be a minimal homeomorphism.
Suppose that $\rh (K_0 (A))$ is dense in $\Aff (T (A)).$
Then $C^* (\Z, X, h)$ is a simple unital \ca\  %
with tracial rank zero
which satisfies the \uct.
\end{thm}

\begin{proof}
It is well known that $C^* (\Z, X, h)$ is simple and unital.

By Theorem~\ref{T1}, for tracial rank zero
it suffices to show that $A_{ \{ y \} }$ has tracial rank zero
for some $y \in X.$
The existence of $y$ with this property
follows from Theorem~\ref{TRZero}.

It follows from Theorem~1.17 of~\cite{RS},
and the description of the class ${\mathcal{N}}$
preceding its statement,
that the crossed product of any commutative \ca\  by any action of $\Z$
satisfies the \uct.
\end{proof}

\begin{cor}\label{C1}
Let $X$ be an infinite compact metric space
with finite covering dimension,
and let $h \colon X \to X$ be a minimal homeomorphism.
Suppose that $\rh (K_0 (A))$ is dense in $\Aff (T (A)).$
Then $C^* (\Z, X, h)$ is a simple \aha\  with no dimension
growth and with real rank zero.
\end{cor}

\begin{proof}
This follows from Theorem~5.2 of~\cite{Lnduke} and Theorem~\ref{TM},
together with the existence of a simple \aha\  %
with the required Elliott invariant.
See Lemma~7.5 of~\cite{Ph11} for how to put everything together.
\end{proof}

The converses of Theorem~\ref{TM} and Corollary~\ref{C1}
are of course easy.
(See Proposition~1.10(a) of~\cite{Ph8}.)

\begin{cor}\label{Tiso}
For $j = 1, 2$ let $X_j$ be an
infinite compact metric space with finite covering dimension,
and let $h_j \colon X \to X$ be a minimal homeomorphism.
Let $A_j = C^*(\Z, X, h_j).$
Suppose that $\rh_{A_j} (K_0 (A_j))$ is dense in $\Aff (T (A_j)))$
for $j = 1, 2.$
Then $A_1 \cong A_2$ if and only
if
\[
(K_0 (A_1), \, K_0 (A_1)_+, \, [1_{A_1}], \, K_1 (A_1))
\cong K_0 (A_2), \, K_0 (A_2)_+, \, [1_{A_2}], \, K_1 (A_2)).
\]
\end{cor}

\begin{proof}
This is immediate
from Theorem~5.2 of~\cite{Lnduke} and Theorem~\ref{TM}.
\end{proof}

\section{Examples}\label{Sec:Examples}

In this section, we describe one general type of example
to which our result applies,
and then prove the isomorphism of four fairly specific pairs
of \tgca s.

For our general example,
we recall briefly the theory of rotation numbers
for automorphisms developed in~\cite{Ex}.
Let $X$ be a connected compact metric space,
let $h \colon X \to X$ be a homeomorphism,
and let $\mu$ be an $h$-invariant Borel measure on $X.$
(Exel works more generally with a unital \ca\  $A,$
an automorphism $\af \colon A \to A,$
and an $\af$-invariant tracial state $\ta$ on $A$ such that
$\ta_* (K_0 (A)) \S \Z.$)
Then the rotation number $\rh_h^{\mu}$
associated with $h$ and $\mu$ (Definition~IV.1 of~\cite{Ex})
is a \hm\  with domain
the kernel of the \hm\  $\id - (h^{-1})^* \colon K^1 (X) \to K^1 (X)$
and codomain $\R / \Z.$
It is defined as follows.
As usual, let $\af \colon C (X) \to C (X)$ be the
automorphism $\af (f) = f \circ h^{-1}.$
Let $u \in U (M_n (C (X))$ satisfy $(\id - (h^{-1})^*) ([u]) = 0.$
Let $v = \af (u^*) u.$
Then $[v] = 0$ in $K_1 (C (X)).$
Increasing the matrix size and replacing $u$ by $\diag (u, 1),$
we may assume that $v \in U_0 (M_n (C (X)).$
Then there exist $a_1, a_2, \ldots, a_m \in M_n (C (X))_{\sa}$
such that $\prod_{k = 1}^m e^{i a_k} = v.$
Now
\[
\rh_h^{\mu} ([u])
 = \Z
   + \frac{1}{2 \pi} \int_X \sum_{k = 1}^m \Tr (a_k (x)) \, d \mu (x).
\]
See Section~IV of~\cite{Ex} for the details, including the proof
that $\rh_h^{\mu}$ is well defined;
see Definition~II.2 of~\cite{Ex} for the determinant used there.

\begin{dfn}\label{IrrRotD}
Let $X$ be a connected compact metric space,
and let $h \colon X \to X$ be a homeomorphism.
We say that $h$ {\emph{has an irrational rotation number}}
if there exists $\et$ in the kernel of
$\id - (h^{-1})^* \colon K^1 (X) \to K^1 (X)$
and an $h$-invariant Borel probability measure $\mu$ on $X$
such that, with $\rh_h^{\mu}$ as above,
one has $\rh_h^{\mu} (\et) \not\in \Q / \Z.$
\end{dfn}

By Theorem~VI.11 of~\cite{Ex},
one only needs to consider $\et$ of the form $[u]$ with $u \in C (X).$

If $X = S^1,$ if $h (\zt) = \exp (2 \pi i \te) \zt$ for $\zt \in S^1,$
and if $\mu$ is normalized Haar measure,
then the unitary $u (\zt) = \zt$
has $\rh_h^{\mu} ([u]) = \te + \Z.$
See Example~IV.5 of~\cite{Ex}.

In this section,
we are particularly interested in the uniquely ergodic case.

\begin{lem}\label{IrrImpDense}
Let $X$ be a connected compact metric space,
let $h \colon  X \to X$ be a minimal homeomorphism,
and let $A = C^* (\Z, X, h).$
Suppose that $h$ is uniquely ergodic
and has an irrational rotation number.
Then $\rh (K_0 (A))$ is dense in $\Aff (T (A)).$
\end{lem}

\begin{proof}
We identify $\Aff (T (A))$ with $\R.$
Clearly $\rh (K_0 (A))$ contains $\Z.$
It follows from Theorem~V.12 of~\cite{Ex} that
the image of $\rh (K_0 (A))$ in $\R / \Z$
contains $\te + \Z$ for some $\te \not\in \Q.$
So $\rh (K_0 (A))$ contains $\Z + \te \Z,$
which is dense in $\R.$
\end{proof}

\begin{cor}\label{IrrImpTR0}
Let $X$ be a finite dimensional connected compact metric space
and let $h \colon X \to X$ be a minimal homeomorphism.
Suppose that $(X, h)$ is uniquely ergodic
and has an irrational rotation number.
Then $C^* (\Z, X, h)$ is a simple unital \aha\  with no
dimension growth and with real rank zero.
\end{cor}

\begin{proof}
Combine Lemma~\ref{IrrImpDense} and Corollary~\ref{C1}.
\end{proof}

\begin{exa}\label{IrratRotAlg}
Corollary~\ref{IrrImpTR0}
allows us to recover the Elliott-Evans Theorem~\cite{EE},
stating that the irrational rotation algebra $A_{\te}$ is
an \ata\  with real rank zero.
As noted after Definition~\ref{IrrRotD},
an irrational rotation has an irrational rotation number,
and it is well known that every irrational rotation is
uniquely ergodic.
So $A_{\te}$ is a simple unital \aha\  with no
dimension growth and with real rank zero.
Since its K-theory is torsion free (easily obtained from
the Pimsner-Voiculescu exact sequence~\cite{PV}),
it follows from known results on classification and the
range of the Elliott invariant that $A_{\te}$ is an \ata.
(See Lemma~7.5 of~\cite{Ph11} for how to put the results together
to get this conclusion.)
\end{exa}

\begin{rmk}\label{IfRR0}
Let $X$ be a connected compact metric space
such that $K^1 (X)$ is finitely generated.
(For example, $X$ could be a finite CW~complex.)
Let $h \colon X \to X$ be a minimal homeomorphism.
Suppose that $C^* (\Z, X, h)$ has tracial rank zero,
or even just suppose that $C^* (\Z, X, h)$ has real rank zero.
Then we claim that $h$ has an irrational rotation number.

To see this, let $A = C^* (\Z, X, h).$
Let $\mu$ be any $h$-invariant Borel probability measure on $X,$
and let $\ta$ be the corresponding trace on $A.$
The map $\rh \colon K_0 (A) \to \Aff (T (A))$ has dense range
by Proposition~1.10(a) of~\cite{Ph8}.
In particular, $\ta_* (K_0 (A))$ is dense in $\R.$
Let $\io \colon C (X) \to A$ be the standard inclusion,
and let $G \S K^1 (X)$ be the kernel of the map
$\id - (h^{-1})^* \colon K^1 (X) \to K^1 (X).$
{}From the Pimsner-Voiculescu exact sequence~\cite{PV},
we get an exact sequence as the top row of the following diagram.
The square commutes by Theorem~V.12 of~\cite{Ex}.
\[
\begin{CD}
0 @>>> \io_* (K_0 (C (X)) @>>> K_0 (A) @>{\partial}>> G @>>> 0 \\
& & & &  @V{\ta_*}VV  @VV{\rh_h^{\mu}}V \\
& & & & \R @>>> \R / \Z.
\end{CD}
\]
Since $X$ is connected, $\ta_* ( \io_* (K_0 (C (X)))) = \Z.$
Therefore $\rh_h^{\mu}$ has dense range.
Since $G$ is finitely generated, this can only happen if its
range is not contained in $\Q / \Z.$
\end{rmk}

We now turn specific examples of isomorphisms of \tgca s.
Except for Example~\ref{IOEx}, these were discussed in~\cite{Ph},
although they predate that paper;
the expected isomorphisms could not be proved there
because of insufficient smoothness.

The first two we consider are covered by Corollary~\ref{IrrImpTR0},
but in fact the computations of the ordered K-theory
have already been done.

\begin{exa}\label{FurstEx}
For each $\te \in [0, 1] \setminus \Q,$
each \cfn\  $f \colon S^1 \to \R,$
and each $d \in \Z \setminus \{ 0 \},$
we define
$h = h_{\te, f, d} \colon S^1 \times S^1 \to S^1 \times S^1$
to be the inverse of the homeomorphism
\[
( \zt_1, \zt_2 ) \mapsto
    {\ts{ \left( \rsz{ e^{2 \pi i \te} \zt_1, \,
  e^{2 \pi i f (\zt_1)} \zt_1^d \zt_2 } \right) }}.
\]
(One sees that the given map does in fact have a
\ct\  inverse by writing down an explicit formula for it.
This homeomorphism is called $\ph_{f, \te}$ in \cite{Kd}.)
The homeomorphism $h$ is minimal by the remark after
Theorem~2.1 in Section~2.3 of \cite{Fr}.
Normalized Lebesgue measure on $S^1 \times S^1$ is invariant,
and when $f$ is Lipschitz this is the only invariant probability measure
(Theorem~2.1 of \cite{Fr}).

Let $B_{\te, f, d} = C^* (\Z, \, S^1 \times S^1, \, h)$
be the crossed product \ca.
The Elliott invariant of $B_{\te, f, d}$ is computed in~\cite{Kd},
and by a more direct method in Example~4.9 of~\cite{Ph7}.
It is described as follows.
Let $A_{\te}$ be the irrational rotation algebra.
The map $S^1 \times S^1 \to S^1$ given by \pj\  on the first
coordinate intertwines $h$ and the irrational rotation
$\zt \mapsto e^{- 2 \pi i \te} \zt,$
and therefore gives a \hm\  $\ph \colon A_{\te} \to B_{\te, f, d}.$
Let $p \in A_{\te}$ be a \pj\  for which the unique tracial state $\ta$
on $A_{\te}$ satisfies $\ta (p) = \te,$
and let $\bt \in K_0 (B_{\te, f, d})$ be the image of the Bott element
in $K_0 (C (S^1 \times S^1))$ under the inclusion.
Then
\[
K_0 (B_{\te, f, d}) \cong \Z [1] \oplus \Z \bt
   \oplus \Z [\ph (p)],
\]
every tracial state $\ta$ on $B_{\te, f, d}$ satisfies
\[
\ta_* (m_1 [1] + m_2 \bt + m_3 [\ph (p)] ) = m_1 + m_3 \te,
\]
the positive part of $K_0$ is given by
\[
K_0 (B_{\te, f, d})_+
 = \{ m_1 [1] + m_2 \bt + m_3 [\ph (p)] \colon
   {\mbox{$m_1 + m_3 \te > 0$ or $m_1 = m_2 = m_3 = 0$}} \},
\]
and
\[
K_1 (B_{\te, f, d}) \cong \Z^3 \oplus \Z / d \Z.
\]
Although this is not explicitly stated in Example~4.9 of~\cite{Ph7},
one sees easily from the computations done there that
the image of $K_1 (C (S^1 \times S^1))$
in $K_1 (B_{\te, f, d})$ is a direct summand,
isomorphic to $\Z \oplus \Z / d \Z,$
with $\Z$ generated by the image of the class of the
unitary $u (\zt_1, \zt_2) = \zt_1$
and $\Z / d Z$ generated by the image of the class of the
unitary $u (\zt_1, \zt_2) = \zt_2.$

Whenever $h$ is uniquely ergodic, Theorem~\ref{T1} applies.
In particular, whenever $f_1$ and $f_2$ are Lipschitz
and $| d_1 | = | d_2 |,$
it follows from Corollary~\ref{Tiso}
that $B_{\te, f_1, d_1} \cong B_{\te, f_2, d_2}.$

When $d = \pm 1,$ the K-theory is torsion free,
so it follows as in Example~\ref{IrratRotAlg} that $B_{\te, f, d}$
is a simple \ata\  with real rank zero.
\end{exa}

\begin{exa}\label{IOEx}
Example~3.6 of~\cite{IO} proves the existence of \mh s
$h_1 \colon X_1 \to X_1$ and $h_2 \colon X_2 \to X_2$
of connected one dimensional compact metric spaces
(which are suspension spaces of \mh s of the Cantor set and
locally look like the product of an interval and the Cantor set)
such that $C^* (\Z, X_1, h_1)$ and $C^* (\Z, X_2, h_2)$
have the same Elliott invariant
but such that $h_1$ and $h_2$ are not flip conjugate.
By Theorem~3.1 and Remark~3.4 of~\cite{BT},
or by Proposition~5.5 of~\cite{LP1},
$h_1$ and $h_2$ are not even topologically orbit equivalent.
However, the Elliott invariants computed in~\cite{IO}
satisfy the dense range condition in Corollary~\ref{Tiso},
so this theorem implies that
$C^* (\Z, X_1, h_1) \cong C^* (\Z, X_2, h_2).$
\end{exa}

The remaining examples involve disconnected spaces,
so arguments using Exel's rotation numbers don't apply.
The first is from Gjerde and Johansen~\cite{GJ},
and is discussed in Example~1.6 of~\cite{Ph}.

\begin{exa}\label{GJEx}
It is shown in Theorem~4 of~\cite{GJ} that every \mh\  $h$ of
the Cantor set $X$ has an almost one to one extension which is a
\mh\  $k$ of a nonhomogeneous space $Y$
with covering dimension $\dim (Y) = 1,$
such that $C^* (\Z, X, h)$ and
$C^* (\Z, Y, k)$ have isomorphic Elliott invariants.
It follows from Corollary~\ref{Tiso} that,
even though $h$ and $k$ are not even topologically orbit equivalent,
we nevertheless have $C^* (\Z, Y, k) \cong C^* (\Z, X, h).$

In particular,
$C^* (\Z, Y, k)$ is a simple \ata\  with real rank zero.
\end{exa}

The following example of Putnam appears as Example~1.7 of~\cite{Ph}.
We refer to~\cite{Ph} for a more detailed description
and the computation of the Elliott invariants of the crossed
products.

\begin{exa}\label{PEx}
For any $\te \in \R \setminus \Q$ let $r_{\te} \colon S^1 \to S^1$
be rotation by $2 \pi \te,$
and let $g_{\te}$ be the \mh\  of a Cantor set $X_{\te}$
obtained by disconnecting $S^1$ along a single orbit of $r_{\te}.$
(This orbit is the unique minimal set of a suitable Denjoy homeomorphism
of $S^1,$ as described in Example~1.7 of~\cite{Ph}.)
Now $\te_1, \, \te_2 \in \R \setminus \Q$ be numbers such that
$1, \, \te_1, \, \te_2$ are linearly independent over $\Q.$
Define \hme s
\[
h_1 = r_{\te_1} \times g_{\te_2}
    \colon S^1 \times X_{\te_2} \to S^1 \times X_{\te_2}
\]
and
\[
h_2 = r_{\te_2} \times g_{\te_1}
 \colon S^1 \times X_{\te_1} \to S^1 \times X_{\te_1}.
\]
Then, as shown in~\cite{Ph}, the \ca s
$C^* (\Z, \, S^1 \times X_{\te_2}, \, h_1)$
and $C^* (\Z, \, S^1 \times X_{\te_1}, \, h_2)$
have the same Elliott invariants,
and they satisfy the dense range condition in Corollary~\ref{Tiso}.
So this theorem implies that
$C^* (\Z, \, S^1 \times X_{\te_2}, \, h_1)
   \cong C^* (\Z, \, S^1 \times X_{\te_1}, \, h_2),$
even though, as shown in Example~1.7 of~\cite{Ph},
$h_1$ and $h_2$ are not topologically orbit equivalent.

Since the K-theory is torsion free (Proposition~1.12 of~\cite{Ph}),
it follows as in Example~\ref{IrratRotAlg}
that the \ca s of such \hme s are
simple \ata s with real rank zero.
\end{exa}

\section{Approximate conjugacy}\label{Sec:AppConj}

Tomiyama has proved (see the corollary at the end of~\cite{Tm})
that two topologically transitive \hme s $h_1$ and $h_2$
of compact metric spaces $X_1$ and $X_2$
are flip conjugate \ifo\  there exists an isomorphism
$\ph \colon C (\Z, X_1, h_1) \to C (\Z, X_2, h_2)$
which sends the canonical copy of $C (X_1)$ in the first
algebra isomorphically to
the canonical copy of $C (X_2)$ in the second algebra.
In particular, this applies to \mh s.
Motivated by this result, the first author introduced
(Definition 3.7 of~\cite{Lnd}) the following definition.

\begin{dfn}\label{DKconj}
Let $h_1 \colon X_1 \to X_1$ and $h_2 \colon X_2 \to X_2$
be minimal homeomorphisms of infinite \cms s $X_1$ and $X_2.$
Suppose that $A_1 = C^* (\Z, X_1, h_1)$ and
$A_2 = C^* (\Z, X_2, h_2)$ have tracial rank zero.
Let $\io_1 \colon C (X_1) \to A_1$ and $\io_2 \colon C (X_2) \to A_2$
be the canonical embeddings.
We say that $h_1$ and
$h_2$ are C*-strongly approximately flip conjugate if there
exists an isomorphism $\ph \colon A_1 \to A_2,$
a homeomorphism $k \colon  X_1 \to X_2,$
and a sequence of unitaries $u_n \in U (A_2)$ such that
\[
\lim_{n \to \infty}
 \| {\mathrm{Ad}} (u_n) \circ \io_1 (f) - \io_2 (f \circ k) \| = 0
\]
for all $f \in C (X_1).$
\end{dfn}

The next theorem is taken from~\cite{Lnd}.
Before stating it, we introduce some terminology.
Recall R{\o}rdam's $KL$-theory,
as defined (for separable \ca s which satisfy the \uct)
after Lemma~5.3 in~\cite{Rr4}.
Suppose that $A_1$ and $A_2$ have tracial rank zero
and satisfy the \uct,
and that $\et \in KL (A_1, A_2)$ gives a unit preserving
order isomorphism $\Gm (\et)_0 \colon K_0 (A_1) \to K_0 (A_2).$
The state space of $K_0 (A_j)$ is exactly the
space of quasitraces on $A_j$ (Theorem 6.9.1 of~\cite{Bl}),
and it is easy to see from the definition
(see 4.9(ix) of~\cite{Lntr1})
that if a \ca\ $A$ has tracial rank zero,
then all quasitraces on $A$ are in fact traces.
So $\et$ induces a unit preserving affine order isomorphism
$L (\et) \colon  \Aff (T (A_1)) \to \Aff (T (A_2)).$

We also define $\rh_A \colon A_{\sa} \to \Aff (T (A))$ by
$\rh_A (a) (\ta) = \ta (a)$ for $\ta \in T (A)$;
no confusion should arise with the map from $K_0 (A)$ with
the same name.

\begin{thm}\label{TKconj}
Let the notation and assumptions be as in Definition~\ref{DKconj}.
Then $h_1$ and $h_2$ are C*-strongly approximately flip conjugate
if and only if there exists a sequence of isomorphisms
$\ch_n \colon C (X_1) \to C (X_2)$
and an element $\et \in KL (A_1, A_2)$ such that:
\begin{itemize}
\item[(1)]
$\Gm (\et) \colon K_* (A) \to K_* (B)$ is a unit
preserving order isomorphism.
\item[(2)]
$\et \times [\io_1] = [\io_2] \times [\ch_n]$
in $KL (C (X_1), A_2)$ for all $n.$
\item[(3)]
$\lim_{n \to \infty}
 \| \rh_{A_2} \circ \io_2 \circ \ch_n (f)
                - L (\et) \circ \io_1 (f) \| = 0$
for all $f \in C (X)_{\sa}.$
\end{itemize}
\end{thm}

\begin{cor}\label{Capp}
Let $h_1 \colon X_1 \to X_1$ and $h_2 \colon X_2 \to X_2$
be minimal homeomorphisms
of finite dimensional infinite \cms s $X_1$ and $X_2.$
Set $A_j = C^* (\Z, X_j, h_j),$
and assume that $\rh_{A_j} (K_0 (A_j))$ is dense in
$\Aff (T (A_j))$ for $j = 1, 2,$
and that $K^* (X_j)$ and $K_* (A_j)$ are torsion free for $j = 1, 2.$
Suppose that there is a unit
preserving order isomorphism
\[
\gm \colon  (K_0 (A_1), K_0 (A_1)_+, [1_{A_1}], K_1 (A_1))
 \to (K_0 (A_2), K_0 (A_2)_+,
[1_{A_2}], K_1 (A_2))
\]
and an isomorphism $\ch \colon C (X_1) \to C (X_2)$
such that
\[
\gm \circ (\io_1)_* = (\io_2)_* \circ \ch_*
\andeqn
L (\gm) \circ \rh_{A_1} \circ \io_1 |_{C (X_1)_{\sa}}
    = \rh_{A_2} \circ \io_2 \circ \ch |_{C (X_1)_{\sa}}.
\]
Then $h_1$ and $h_2$ are C*-strongly approximately flip conjugate.
\end{cor}

\begin{proof}
By Theorem~\ref{TM},
the algebras $A_1$ and $A_2$ have tracial rank zero.
Since $K_* (A_j)$ is torsion free, the \uct\  implies that
$KL (A_1, A_2) = \Hom (K_* (A_1), \, K_* (A_2)).$
(See~the discussion before Lemma~5.2 and
before and after Lemma~5.3 in~\cite{Rr4}.)
Similarly, since $K_* (C (X_1))$ is torsion free,
$KL (C (X_1), A_2) = \Hom (K_* (C (X_1), \, K_* (A_2)).$
Now the corollary follows from Theorem~\ref{TKconj}.
\end{proof}

\begin{exa}\label{VEfurst2}
Let $\te \in [0, 1] \setminus \Q,$
let $f_1, f_2 \colon S^1 \to \R$ be \cfn s,
and let $d_1, d_2 \in \Z \setminus \{ 0 \}$
satisfy  $|d_1| = |d_2|.$
Let $h_1 = h_{\te, f_1, d_1}$ and $h_2 = h_{\te, f_2, d_2}$
be uniquely ergodic Furstenberg transformations on $X = S^1 \times S^1$
as in Example~\ref{FurstEx},
and let $A_1 = C^* (\Z, X, h_1)$ and $A_2 = C^* (\Z, X, h_2).$
Let $\io_1 \colon  C (X) \to A_1$ and
$\io_2 \colon  C (X) \to A_2$ be the canonical embeddings.
{}From the K-theory computation described in Example~\ref{FurstEx},
one sees that there is an isomorphism
$\gm \colon K_* (A_1) \to K_* (A_2)$
which is an order isomorphism on $K_0$ and such that the following
diagram commutes:
\[
\begin{CD}
K_* (C (X)) @>{=}>> K_* (C (X)) \\
 @V{(\io_1)_*}VV  @VV{(\io_2)_*}V \\
K_* (A_1) @>>{\gm}> K_* (A_2).
\end{CD}
\]
It follows from Theorem~5.2 of~\cite{Lnduke} that there is an isomorphism
$\ph \colon A_1 \to A_2$ such that $\ph_* = \gm.$

Let $\ta_1$ and $\ta_2$ be the unique tracial states on $A_1$ and $A_2.$
We claim that $\ta_2 \circ \ph \circ \io_1 = \ta_2 \circ \io_2.$
By uniqueness, we have $\ta_2 \circ \ph = \ta_1.$
Let $\mu$ be normalized Lebesgue measure on $S^1 \times S^1,$
which is the unique ergodic measure for both $h_2$ and $h_2.$
Then, for $f \in C (X),$ we have
\[
\ta_2 \circ \ph (\io_1 (f))
 = \ta_1 ( \io_1 (f))
 = \int_{S^1 \times S^1} f \, d \mu
 = \ta_2 (\io_2 (f)).
\]
This proves the claim.
It now follows from Corollary~\ref{Capp} that $h_1$ and $h_2$ are
C*-strongly approximately flip conjugate.

Suppose that $f_1 = 0$ and $d_1 = d_2 = 1.$
It has been shown in Example~2.1 of~\cite{Ph}
that there are choices of $f_2$ for which
$h_1$ and $h_2$ are not flip conjugate,
and even not topologically orbit equivalent.
\end{exa}

\begin{rmk}
Let $X$ be the Cantor set and $h_1, h_2 \colon X \to X$ be two
minimal homeomorphisms.
Then Corollary~\ref{Capp}
implies that $h_1$ and $h_2$ are C*-strongly
approximately flip conjugate if and only if
there is a unit preserving order isomorphism between
the K-groups $K^0 (X, h_1)$ and $K^0 (X, h_2)$
of Definition~1.11 of~\cite{GPS}.
In~\cite{LM1}, it is shown this is also equivalent
to the existence of sequences
$(k_n)_{n \in \N}$ and $(l_n)_{n \in \N}$
of homeomorphisms from $X$ to $X$ such that:
\begin{itemize}
\item[(1)]
For all $f \in C (X),$ we have
\[
\lim_{n \to \infty}
  \|f \circ k_n^{-1} \circ h_2 \circ k_n - f \circ h_1 \| = 0
\andeqn
\lim_{n \to \infty}
  \|f \circ l_n^{-1} \circ h_1 \circ l_n - f \circ h_2 \| = 0
\]
\item[(2)]
For every \pj\  $p \in C (X),$
for all sufficiently large $n,$
the images of $k_n^* ([p])$ and $k_{n+1}^* ([p])$
in $K_0 (C^* (\Z, X, h_1)$ are equal;
similarly,
for every \pj\  $q \in C (X),$
for all sufficiently large $n,$
the images of $l_n^* ([q])$ and $l_{n+1}^* ([q])$
in $K_0 (C^* (\Z, X, h_2)$ are equal.
\end{itemize}
If such $(k_n)_{n \in \N}$ and $(l_n)_{n \in \N}$ exist,
we say $(X, h_1)$ and $(X, h_2)$ are approximately $K$-conjugate.
Approximate $K$-conjugacy is also equivalent
to strong orbit equivalence as in~\cite{GPS}.
(See~\cite{LM2}.)
\end{rmk}

Approximate $K$-conjugacy can also be defined for more general spaces.
See~\cite{LM1}, \cite{Lnd}, and~\cite{LM2}.
In~\cite{LM2}, approximate $K$-conjugacy has been studied for
\hme s of the product of the Cantor set and the circle.

We would like to ask the following question.
Let $(X_1, h_1)$ and $(X_2, h_2)$ be two minimal
systems,
and let $A_1 = C^* (\Z, X_1, h_1)$ and $A_2 = C^* (\Z, X_2, h_2).$
Suppose that $\rh (K_0 (A_1))$ and $\rh (K_0 (A_2))$ are dense
in $\Aff (T (A_1))$ and $\Aff (T (A_2)),$
and that $h_1$ and $h_2$ are
C*-strongly approximately flip conjugate.
What additional hypotheses are required for
$h_1$ and $h_2$ to be approximately $K$-conjugate?

\flushleft{\textit{\mbox{} \\
Huaxin Lin\\
email: hxlin@noether.uoregon.edu\\
Department of Mathematics\\
East China Normal University\\
Shanghai, China\\
and (current)\\
Department of Mathematics\\
University of Oregon\\
Eugene, Oregon 97403-1222\\
U.S.A.\\
\mbox{} \\
N.\  Christopher Phillips\\
email: ncp@darkwing.uoregon.edu\\
Department of Mathematics\\
University of Oregon\\
Eugene, Oregon 97403-1222\\
U.S.A.\\
}}


\begin{thebibliography}{99}

\bibitem{Bg} I.\  D.\  Berg,
{\emph{On approximation of normal operators by weighted shifts}},
Michigan Math.\  J.\  {\textbf{21}}(1974), 377--383.

\bibitem{Bl} B.\  Blackadar,
{\emph{K-Theory for Operator Algebras}},
MSRI Publication Series {\textbf{5}}, Springer-Verlag,
New York, Heidelberg, Berlin, Tokyo, 1986.

\bibitem{BT} M.\  Boyle and J.\  Tomiyama,
{\emph{Bounded topological orbit equivalence and C*-algebras}},
J.\  Math.\  Soc.\  Japan {\textbf{50}}(1998), 317--329.

\bibitem{Br} L.\  G.\  Brown,
{\emph{Interpolation by projections in C*-algebras of real rank zero}},
J.\  Operator Theory {\textbf{26}}(1991), 383--387.

\bibitem{Cn0} A.\  Connes,
{\emph{An analogue of the Thom isomorphism for crossed products of
a C*-algebra by an action of ${\mathbb{R}}$}},
Advances in Math.\  {\textbf{39}}(1981), 31--55.

\bibitem{El5} G.\  A.\  Elliott,
{\emph{The classification problem for amenable C*-algebras}},
pages 922--932 in:
{\emph{Proceedings of the International Congress of Mathematicians,
Z\"{u}rich, 1994}}, S.\  D.\  Chatterji, ed.,
Birkh\"{a}user, Basel, 1995.

\bibitem{EE} G.\  A.\  Elliott and D.\  E.\  Evans,
{\emph{The structure of the irrational rotation algebra}},
Ann.\  of Math.\  (2) {\textbf{138}}(1993), 477--501.

\bibitem{EG} G.\  A.\  Elliott and G.\  Gong,
{\emph{On the classification of C*-algebras of real rank zero, II}},
Ann.\  Math.\  {\textbf{144}}(1996), 497--610.

\bibitem{Ex} R.\  Exel,
{\emph{Rotation numbers for automorphisms of C*-algebras}},
Pacific J.\  Math.\  {\textbf{127}}(1987),  31--89.

\bibitem{Fr} H.\  Furstenberg,
{\emph{Strict ergodicity and transformation of the torus}},
Amer.\  J.\  Math.\  {\textbf{83}}(1961), 573--601.

\bibitem{GJ} R.\  Gjerde and {\O}.\  Johansen,
{\emph{C*-algebras associated to non-homogeneous minimal systems and
their K-theory}}, Math.\  Scand.\  {\textbf{85}}(1999), 87--104.

\bibitem{GPS}
T.\  Giordano, I.\  F.\  Putnam, and C.\  F.\  Skau,
{\emph{Topological orbit equivalence and C*-crossed products}},
J.\  reine angew.\  Math.\  {\textbf{469}}(1995), 51--111.

\bibitem{HPS}
R.\  H.\  Herman, I.\  F.\  Putnam, and C.\  F.\  Skau,
{\emph{Ordered Bratteli diagrams, dimension groups, and
topological dynamics}},
International J.\  Math.\  {\textbf{3}}(1992), 827--864.

\bibitem{HLX} S.\  Hu, H.\  Lin, and Y.\  Xue,
{\emph{Tracial topological rank of extensions of C*-algebras}},
Math.\  Scand.\  {\textbf{94}}(2004), 125--147.

\bibitem{IO} B.\  Itz\'{a}-Ortiz,
{\emph{The C*-algebras associated to time-$t$ automorphisms of
mapping tori}},
in preparation.

\bibitem{Ji} R.\  Ji,
{\emph{On the Crossed Product C*-Algebras Associated with
Furstenberg Transformations on Tori}},
Ph.D.\  Thesis, State University of New York at Stony Brook, 1986.

\bibitem{Kd} K.\  Kodaka,
{\emph{The positive cones of $K_0$-groups of crossed products
associated with Furstenberg transformations on the $2$-torus}},
Proc.\  Edinburgh  Math.\  Soc.\  {\textbf{43}}(2000), 167--175.

\bibitem{Kl} J.\  Kulesza,
{\emph{Zero-dimensional covers of finite dimensional
dynamical systems}},
Ergod.\  Th.\  Dynam.\  Sys.\  {\textbf{15}}(1995), 939--950.

\bibitem{Lnnaf} H.\  Lin,
{\emph{Approximation by normal elements with
finite spectra in simple AF-algebras}},
J.\  Operator Theory {\textbf{31}}(1994), 83--98.

\bibitem{Ln1} H.\  Lin, {\emph{Tracially AF C*-algebras}},
Trans.\  Amer.\  Math.\  Soc.\  {\textbf{353}}(2001), 693--722.

\bibitem{Ln2} H.\  Lin,
{\emph{The tracial topological rank of C*-algebras}},
Proc.\  London\  Math.\  Soc.\  {\textbf{83}}(2001),  199--234.

\bibitem{Lnduke} H.\  Lin,
{\emph{Classification of simple C*-algebras with
tracial topological rank zero}},
Duke Math.\  J., to appear.

\bibitem{Lntr1} H.\ Lin, {\emph{Simple nuclear C*-algebras of tracial 
topological rank one}}, preprint.

\bibitem{Lnd} H.\  Lin,
{\emph{Classification of homomorphisms and minimal dynamical systems}},
preprint.

\bibitem{LM1} H.\  Lin and H.\  Matui,
{\emph{Minimal dynamical systems and approximate conjugacy}},
preprint (math.OA/0402309).

\bibitem{LM2}  H.\  Lin and H.\  Matui,
{\emph{Minimal dynamical systems on the product of the Cantor
set and the circle}}, preprint.

\bibitem{LO} H.\  Lin and H.\  Osaka,
{\emph{The Rokhlin property and the tracial topological rank}},
J.\  Funct.\  Anal., to appear (math.OA/0402094).

\bibitem{Lnq} Q.\  Lin,
{\emph{Analytic structure of the transformation group C*-algebra
associated with minimal dynamical systems}},
preprint.

\bibitem{LP1} Q.\  Lin and N.\  C.\  Phillips,
{\emph{Ordered K-theory for C*-algebras of minimal homeomorphisms}},
pages 289--314 in: {\emph{Operator Algebras and Operator Theory}},
L.\  Ge, etc.\ (eds.),
Contemporary Mathematics vol.\  228, 1998.

\bibitem{LP2} Q.\  Lin and N.\  C.\  Phillips,
{\emph{Direct limit decomposition for C*-algebras
of minimal diffeomorphisms}},
pages 107--133 in:
{\emph{Operator Algebras and Applications,}}
Advanced Studies in Pure Mathematics vol.\  38,
Mathematical Society of Japan, 2004.

\bibitem{LP3} Q.\  Lin and N.\  C.\  Phillips,
{\emph{The structure of C*-algebras of minimal diffeomorphisms}},
in preparation.

\bibitem{Pc86} J.\  A.\  Packer,
{\emph{K-theoretic invariants for C*-algebras associated to
  transformations and induced flows}},
J.\  Funct.\  Anal.\  {\textbf{67}}(1986), 25--59.

\bibitem{Ph} N.\  C.\  Phillips,
{\emph{When are crossed products by minimal diffeomorphisms
isomorphic?}}, pages 341--364 in:
{\emph{Operator Algebras and Mathematical Physics}}
(Conference Proceedings, Constan\c{t}a, (Romania) July 2--7, 2001),
J.-M.\  Combes, J.\  Cuntz, G.\  A.\  Elliott, G.\  Nenciu,
H.\  Seidentop, \c{S}.\  Str\v{a}til\v{a} (eds.),
The Theta Foundation, Bucharest, 2003.

\bibitem{Ph6}  N.\  C.\  Phillips,
{\emph{Recursive subhomogeneous algebras}}, preprint (math.OA/0101156).

\bibitem{Ph7}  N.\  C.\  Phillips,
{\emph{Cancellation and stable rank for direct limits of recursive
subhomogeneous algebras}}, preprint (math.OA/0101157).

\bibitem{Ph8}  N.\  C.\  Phillips,
{\emph{Real rank and property~(SP) for direct limits of recursive
subhomogeneous algebras}}, preprint.

\bibitem{Ph11}  N.\  C.\  Phillips,
{\emph{Crossed products by finite cyclic group actions
  with the tracial Rokhlin property}},
preprint.

\bibitem{PV} M.\  Pimsner and D.\  Voiculescu,
{\emph{Exact sequences
for K-groups and Ext-groups of certain cross-products of C*-algebras}},
J.\  Operator Theory {\textbf{4}}(1980), 93--118.

\bibitem{Pp} S.\  Popa,
{\emph{On local finite dimensional approximation of C*-algebras}},
Pacific J.\  Math.\  {\textbf{181}}(1997), 141--158.

\bibitem{Pt1} I.\   F.\  Putnam,
{\emph{The C*-algebras associated with minimal
homeomorphisms of the Cantor set}},
Pacific J.\  Math.\  {\textbf{136}}(1989), 329--353.

\bibitem{Pt2} I.\   F.\  Putnam,
{\emph{On the topological stable rank of certain
transformation group C*-algebras}},
Ergod.\  Th.\  Dynam.\  Sys.\  {\textbf{10}}(1990), 197--207.

\bibitem{RM} K.\  Reihani and P.\  Milnes,
{\emph{C*-algebras from Anzai flows and their K-groups}},
preprint.

\bibitem{Rr4} M.\  R{\o}rdam,
{\emph{Classification of certain infinite simple C*-algebras}},
J.\  Funct.\  Anal.\  {\textbf{131}}(1995),  415--458.

\bibitem{RS}  J.\  Rosenberg and C.\  Schochet,
{\emph{The K\"{u}nneth theorem and the universal coefficient theorem
for Kasparov's generalized K-functor}},
Duke Math.\  J.\  {\textbf{55}}(1987), 431--474.

\bibitem{Tm} J.\  Tomiyama, {\emph{Topological full groups
and structure of normalizers in transformation group C*-algebras}},
Pacific J.\  Math.\  {\textbf{173}}(1996), 571--583.

\end{thebibliography}
\end{document}